\documentclass{elsarticle}
 
\usepackage{amsmath,amsfonts,amssymb,amsthm}
\usepackage{stmaryrd,mathrsfs}
\usepackage{hyperref}
\usepackage[all]{hypcap}
\usepackage[font=small,labelsep=colon]{caption}
\usepackage{rotating}

\newtheorem{thm}{Theorem}[section]
\newtheorem{lem}[thm]{Lemma}
\newtheorem{cor}[thm]{Corollary}
\newtheorem{prop}[thm]{Proposition}
\newtheorem{bigthm}{Theorem}
\newtheorem*{thmApreview}{Theorem \ref{thm: delta det by 0s}}
\newtheorem*{thmBpreview}{Theorem \ref{max-corner}}
\newtheorem*{thmCpreview}{Theorem \ref{bounds1}}
\newtheorem{thmapp}{Theorem}
\newtheorem{lemapp}[thmapp]{Lemma}

\newdefinition{defn}[thm]{Definition}
\newdefinition{rem}[thm]{Remark}
\newdefinition{exmp}[thm]{Example}
\newdefinition{remapp}[thmapp]{Remark}
\newdefinition{defnapp}[thmapp]{Definition}

\newproof{pf}{Proof}
\newproof{potC}{Proof of Theorem \ref{bounds1}}
\newproof{potC2}{Proof of Theorem \ref{bounds1} (for $p=2$)}
\newproof{potC3}{Proof of Theorem \ref{bounds1} (for special cell classes)}
\newproof{proofHKMformula}{Proof of Theorem \ref{HKM versus delta}}
\newproof{potLipschitz}{Proof of Theorem \ref{mu is lipschitz}}

\DeclareMathOperator{\charac}{char}
\DeclareMathOperator{\genus}{genus}
\DeclareMathOperator{\syz}{Syz}
\DeclareMathOperator{\proj}{Proj}
\DeclareMathOperator{\slope}{slope}
\DeclareMathOperator{\length}{length}

\newcommand{\ideal}[1]{{\langle{#1}\rangle}}
\newcommand{\cclass}[1]{{\left[{#1}\right]}}
\newcommand{\cell}[1]{{\left({#1}\right)}}
\newcommand{\vv}[1]{{\mathbf{{#1}}}}
\newcommand{\norm}[1]{{\|\vv{{#1}}\|}}
\newcommand{\vb}[1]{{\mathcal{#1}}}

\newcommand{\Xset}{{\mathcal{X}}}
\newcommand{\ZEROset}{{\mathcal{Z}}}
\newcommand{\projcurve}{{Y}}
\newcommand{\frob}{{\mathcal{F}}}
\newcommand{\diag}{{\Delta_1}}
\newcommand{\slice}{{\Delta_2}}
\newcommand{\mm}{{\mathfrak{m}}}

\newcommand{\eg}{\emph{e.g.}, }
\newcommand{\ie}{\emph{i.e.}, }

\begin{document}

\begin{frontmatter}
\title{Syzygy gap fractals---I. \\ Some structural results and an upper bound}
\author[pedro]{Pedro Teixeira}
\address[pedro]{Knox College,  2 E. South Street, Galesburg, IL  61401-4999, USA}
\journal{Journal of Algebra}
\ead{pteixeir@knox.edu}

\begin{abstract}
 $\Bbbk$ is a field of characteristic $p>0$, and $\ell_1,\ldots,\ell_n$ are  linear forms in $\Bbbk[x,y]$.
Intending applications to   Hilbert--Kunz theory,   to each triple  $C=(F,G,H)$ of nonzero homogeneous  elements of  $\Bbbk[x,y]$ we associate a function $\delta_C$  that encodes the ``syzygy gaps''  of $F^q$, $G^q$, and $H^q\ell_1^{a_1}\cdots \ell_n^{a_n}$, for all $q=p^e$ and $a_i\le q$. 
These are close relatives  of  functions introduced in  \emph{ $p$-Fractals and power series---I} [P. Monsky, P. Teixeira,  $p$-Fractals and power series---I.
Some 2 variable results, J. Algebra  280 (2004) 505--536].  Like  their relatives, the $\delta_C$    exhibit surprising  self-similarity   related to ``magnification by   $p$,'' and knowledge of their structure allows the explicit computation of various Hilbert--Kunz functions.   
 
We   show that these ``syzygy gap fractals'' are  determined by  their zeros and  have a simple  behavior near their local maxima, and  derive an  upper bound for their  local maxima  which has long been  conjectured by Monsky.  Our results will allow us, in a sequel to this  paper,  to determine the structure of the $\delta_C$ by studying the vanishing of certain determinants.   \end{abstract}

\end{frontmatter}

\section{Introduction}

Let $\Bbbk$ be a field of characteristic $p>0$ and $A=\Bbbk[x,y]$. Let $F$, $G$, and $H\in A$ be nonzero homogeneous polynomials  with no common factor. The module of syzygies of $F$, $G$, and $H$ is   free   on two homogeneous generators; let  $\alpha\ge \beta$ be their degrees. We define $\delta(F,G,H)=\alpha-\beta$; this  is   the \emph{syzygy gap}  of $F$, $G$, and $H$. Syzygy gaps were introduced by Han  \cite{han thesis}, were    studied by the author in his thesis \cite{tese}, and have since made   scattered appearances in the literature    \cite{brenner and kaid,hara, mason}. 

This  paper is concerned with a family of   functions  introduced in \cite{tese}, defined in terms of syzygy gaps.  Fix pairwise prime  linear forms $\ell_1,\ldots,\ell_n\in A$, and let $C=(F,G,H)$ be a triple of nonzero homogeneous elements of $A$ such that $F$, $G$, and $H\ell_1\cdots\ell_n$ have no common factor. Let  $\mathscr{I}=[0,1]\cap\mathbb{Z}[1/p]$.  We define   $\delta_C:\mathscr{I}^n\to \mathbb{Q}$  as follows:
for  each  $q=p^e$ and    $\vv{a}=(a_1,\ldots,a_n)\in \mathbb{Z}^n$ with $0\le a_i\le q$, we set  
\[\delta_C\left(\frac{\vv{a}}{q}\right)=\frac{1}{q}\cdot\delta(F^q,G^q,H^q\ell_1^{a_1}\cdots \ell_n^{a_n}).\]

\begin{figure}
\centering
\includegraphics[height=3.5in]{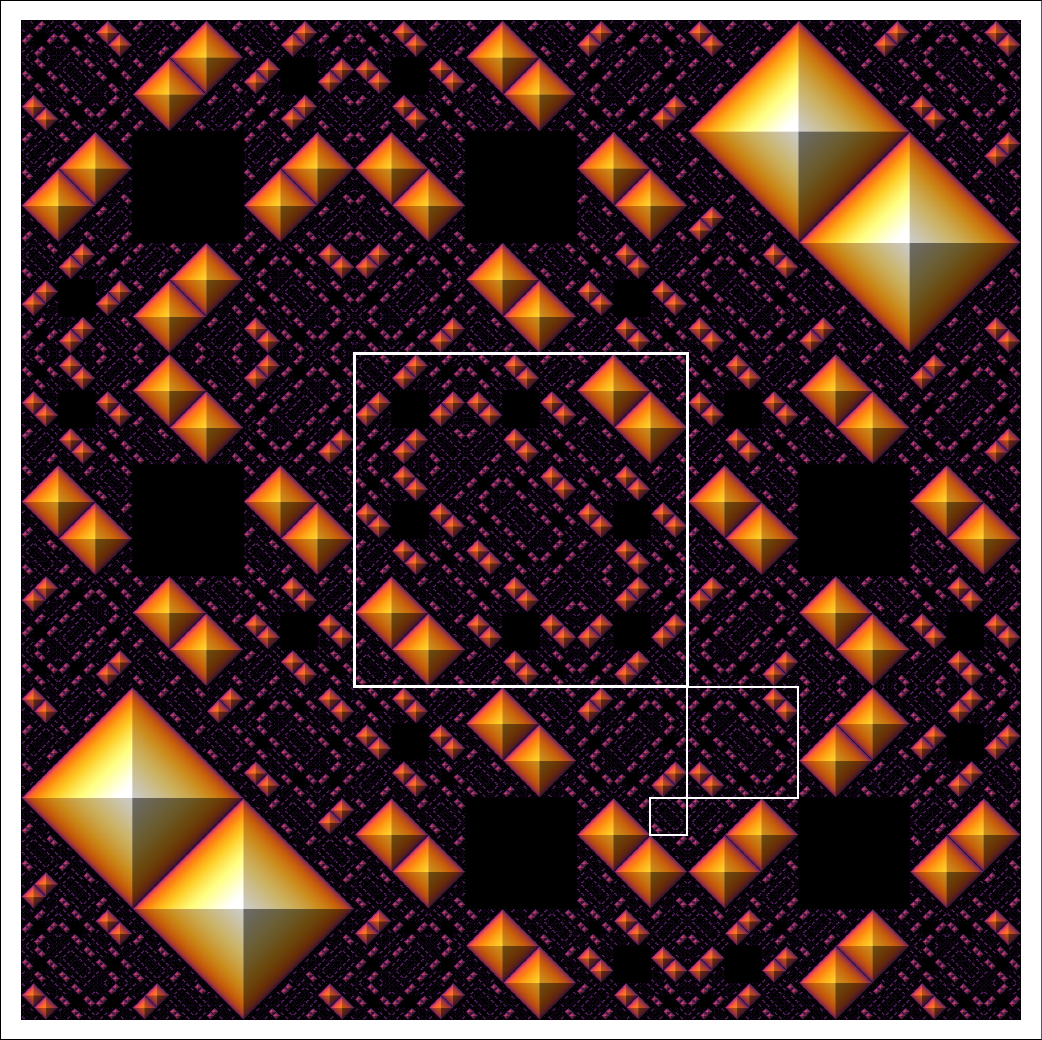} 
\caption{A syzygy gap fractal in characteristic 3}
\label{fig: alpha}
\end{figure}

A two-dimensional ``slice'' of one such function is shown in Figure \ref{fig: alpha},  as a relief plot---zeros are shown in black, and  other   values are encoded by  color (higher value $\leftrightarrow$ lighter color).
 The white  squares  highlight three smaller copies of   the plot    contained within itself.  The   $\delta_C$ often bear  this kind of self-similarity, and if $\Bbbk$ is finite they are  \emph{$p$-fractals}, in the sense of  \cite{pfractals1}. As such,  they can be   characterized   by a finite set of values and a finite set of functional equations---the ``magnification rules''---that prescribe  how    pieces patch together to form the function.

 The ``syzygy gap fractals'' $\delta_C$  are closely related to  functions $\varphi_I:\mathscr{I}^n\to \mathbb{Q}$ introduced  (in  a more general setting) by Monsky and the author  in \cite{pfractals1}. Restricting to the situation at hand  we let $I=\ideal{F,G}:H$ and define, for $\vv{a}$ and $q$ as above, 
\[\varphi_I\left(\frac{\vv{a}}{q}\right)=\frac{1}{q^2}\cdot\deg \ideal{I^{[q]},\ell_1^{a_1}\cdots \ell_n^{a_n}},\]
 where $I^{[q]}=\ideal{u^q:u\in I}$ and $\deg$ denotes the degree or colength of an ideal. 
Then $\delta_C^2$ and $4\varphi_I$ differ by a   polynomial in the coordinate functions (see Eq. (\ref{phi_I and delta_C}) in Section \ref{gaps-and-fractals}). 

In \cite{pfractals2} the  $\varphi_I$  are used in the proof of rationality and computation of the Hilbert--Kunz series and multiplicities of power series of the form $f_1(x_1,y_1)+\cdots+f_m(x_m,y_m)$ with coefficients in a finite field. More specifically, the  ``$p$-fractalness'' of the $\varphi_I$,  established in \cite{pfractals1}, gives  us the rationality result, while   knowledge of the  magnification rules  for those functions  (when available)  allows us to explicitly compute related  Hilbert--Kunz series and multiplicities.  In the present paper we focus on the homogeneous case and find   properties of the syzygy gap fractals $\delta_C$ that will help in those explicit calculations. 

 The main results of this paper concern  the zeros and local maxima of the  $\delta_C$. We  prove that   these functions (when nontrivial) are determined by their zeros:  

\begin{thmApreview} 
 Let $\ZEROset=\{\vv{z}\in \mathscr{I}^{n}\,|\, \delta_C(\vv{z})=0\}$. 
\begin{itemize}
\item If $\ZEROset$ is empty, then   $\delta_C$ is linear; it takes on a minimum value at a corner  $\vv{u}$ of $\mathscr{I}^n$ and,  at each     $\vv{t}\in \mathscr{I}^n$, \[\delta_C(\vv{t})=\delta_C(\vv{u})+d(\vv{t},\vv{u}),\]
 where $d(\vv{t},\vv{u})$ is the   taxi-cab distance between $\vv{t}$ and $\vv{u}$.
\item
If $\ZEROset$ is nonempty then $\delta_C(\vv{t})$ is the  taxi-cab distance from
$\vv{t}$ to the set $\ZEROset$, for all $\vv{t}\in \mathscr{I}^{n}$.
\end{itemize}
\end{thmApreview}

This result has some interesting consequences that will be explored in a sequel to this paper: since the vanishing of  the syzygy gap is tied to the vanishing of a certain determinant in the coefficients of the polynomials, we shall use those determinants in the investigation of the $\delta_C$. We shall prove that the $\delta_C$ are completely determined by \emph{finitely many} such determinants, and this will give us a powerful tool for determining magnification rules, thus allowing  the explicit (and even automatic) calculation of various Hilbert--Kunz series and multiplicities. 
 
Related to Theorem \ref{thm: delta det by 0s} is our next   result, which shows  that  each local maximum of   $\delta_C$ determines the  behavior of the function on a certain neighborhood:

\begin{thmBpreview} Let  $q$ be a power of $p$, and let ${\Xset_{q}}$ be the set consisting of all points of $\mathscr{I}^n$ with denominator $q$.  Suppose the restriction of  $\delta_C$ to $\Xset_q$ attains 
a ``local maximum'' at  $\vv{u}$, in the sense that the values 
of $\delta_C$ at all  points of $\Xset_q$ adjacent to $\vv{u}$ are  smaller  than $\delta_C(\vv{u})$.  Then 
\[\delta_C(\vv{t})=\delta_C(\vv{u})-d(\vv{t},\vv{u}),\]
for all $\vv{t}\in \mathscr{I}^{n}$ with $d(\vv{t},\vv{u})\le \delta_C(\vv{u})$.  In particular, $\delta_C$ is piecewise linear on that region, and has a local maximum at $\vv{u}$ in the usual sense. 
\end{thmBpreview}

Theorem \ref{max-corner} plays a major role in understanding the structure of the $\delta_C$, and is fundamental in the proof  of the last of our results, which shows the existence of a certain upper bound for the  $\delta_C$ at  their  local maxima:

\begin{thmCpreview} 
Suppose $\delta_C$ has a local maximum at $\vv{a}/q$, where $q>1$ and  some $a_i$ is not divisible by $p$. Then    
\[\delta_C\left(\frac{\vv{a}}{q}\right)\le \frac{n-2}{q}.\] 
\end{thmCpreview}
 
 This bound has long been conjectured by Monsky; in \cite{mason} he proved it holds when $C=(x, y,1)$.  The approach used here follows   closely an alternate, unpublished  proof by Monsky of his result from \cite{mason}, where  he gets  information on the local maxima of $\delta_C$  by combining a theorem of Trivedi \cite[Theorem 5.3]{trivedi}  on the Hilbert--Kunz multiplicity of a certain projective plane curve and a formula expressing that  same  multiplicity in terms of a continuous extension of $\delta_C$. Here we combine results of Brenner \cite[Corollary 4.4]{brenner} and Trivedi \cite[Lemma 5.2]{trivedi} and follow essentially the same track  to get to the stronger result, modulo some technical obstacles.

This paper is structured as follows. In Section \ref{gaps-and-delta} we prove some    properties of syzygy gaps independent of the characteristic.  Starting in Section \ref{gaps-and-fractals} we restrict our attention to  positive characteristic; we introduce  the functions $\delta_C$ and look at various examples, and in Section \ref{zeros} we  prove  Theorems \ref{thm: delta det by 0s} and \ref{max-corner}. In Section \ref{operators} we introduce     operators on the ``cells'' $C=(F,G,H)$  that  are mirrored by ``magnifications'' and ``reflections'' on the corresponding functions. While the   ``$p$-fractalness'' of the $\delta_C$ when  $\Bbbk$ is finite is not directly relevant to this paper,  it follows without much effort from the machinery introduced in Section \ref{operators}, so we   present a  proof   in that section. Finally, in Section \ref{bounds} we prove Theorem \ref{bounds1}.

Throughout this paper $p$ denotes a prime number and (lower-case)   $q$ is  used   exclusively for powers of $p$;   $\Bbbk$ is a field, assumed everywhere but in Section \ref{gaps-and-delta}  to be of characteristic $p$;  $\mathscr{I}$ is the set of rational numbers in $[0,1]$ whose denominators are powers of $p$.

\section{Syzygy gaps}\label{gaps-and-delta}

Throughout this  section $\Bbbk$ is a field of arbitrary characteristic, and  $F$, $G$, and $H$  are nonzero homogeneous elements of $A=\Bbbk[x,y]$.  
By the Hilbert Syzygy Theorem, the module of syzygies  of  $(F,G, H)$, denoted by $\syz(F,G,H)$, is free on two homogeneous generators.

  \begin{defn}\label{defsyzgap} 
The \emph{syzygy gap} of $F$, $G$, and $H$
is the nonnegative integer $\delta=n-m$, where $m\le n$ are the degrees 
of the generators of $\syz(F,G,H)$. 
\end{defn}

In  this section we   prove some general properties of syzygy gaps that are characteristic independent. Some  of these  appeared  in \cite{mason}, but are included here, with proofs, for  completeness. Our first result  relates the syzygy gap to the degree of the ideal $\ideal{F,G,H}$ when this degree is finite. 
 
\begin{prop}\label{formula1}Let $F$, $G$, and $H\in \Bbbk[x,y]$ be nonzero homogeneous polynomials with no common factor, of degrees $d_1$, $d_2$, and $d_3$, and let $\delta$ be their syzygy gap.  Then 
\[4\deg\ideal{F,G,H}  =Q(d_1,d_2,d_3) +\delta^2,\] 
where 
\[Q(d_1,d_2,d_3)  =
2(d_1d_2+d_1d_3+d_2d_3)-{d_1}^2-{d_2}^2-{d_3}^2.\] 
\end{prop}

\begin{pf} Let $m$ and $n$ be as in Definition \ref{defsyzgap}.
$A/\ideal{F,G,H}$ has a graded  free resolution
\[0\rightarrow A(-m)\oplus A(-n)\rightarrow \bigoplus_{i=1}^3A(-d_i)\rightarrow A\rightarrow A/\ideal{F,G,H}\rightarrow 0, \]
so  the Hilbert series of $A/\ideal{F,G,H}$ is
\[h(t)=\frac{1-t^{d_1}-t^{d_2}-t^{d_3}+t^{m}+t^{n}}{(1-t)^2}.\]
Since $F$, $G$, and $H$ have no common factor, $h(1)=\deg\ideal{F,G,H}$ is finite. 
Differentiating $(1-t)^2h(t)$ and setting $t=1$ we find 
\begin{equation}\label{relation degrees vs pol degrees}
m+n=d_1+d_2+d_3.
\end{equation} 
Differentiating  $(1-t)^2h(t)$  twice and setting $t=1$ we get
\[2h(1)=-d_1(d_1-1)-d_2(d_2-1)-d_3(d_3-1)+m(m-1)+n(n-1),\] 
and the result follows easily.
\qed
\end{pf}

Eq. (\ref{relation degrees vs pol degrees}) shows that $\delta=d_1+d_2+d_3-2m$; this suggests the following definition:

\begin{defn} 
 $\delta(F,G,H)=\deg F+\deg G+\deg H-2m(F,G,H),$ where $m(F,G,H)$ is the least degree of a nontrivial syzygy of  $(F,G,H)$.  
\end{defn}

\begin{rem} If $F$, $G$, and $H$ have no common factor, $\delta(F,G,H)$ is just the syzygy gap of $F$, $G$, and $H$. In this case, Proposition \ref{formula1} shows  that $\delta(F,G,H)$ remains unchanged under any modification  in the polynomials $F$, $G$, and $H$ that fixes their degrees  and  the ideal $\ideal{F,G,H}$ or, more generally, that fixes $Q(\deg F,\deg G,\deg H)$ and $\deg\ideal{F,G,H}$. 
\end{rem}

\begin{rem}\label{triangle-inequalities}
If $d_{3}\ge d_{1} +d_{2}$ and $F$ and  $G$ have no common factor,  then   $(G,-F,0)$ is a syzygy   of minimal degree, and 
$\delta(F,G,H)=d_{3}-d_{1}-d_{2}.$
\end{rem}
 
 \begin{prop}\label{deltaproperties1} Let $P\in \Bbbk[x,y]$ be a nonzero homogeneous polynomial. Then
 \begin{enumerate}
  \item $\delta(PF,PG,PG)=\delta(F,G,H)+\deg P$;
  \item  $\delta(PF,PG,H)=\delta(F,G,H)$, whenever    $P$ is  prime to $H$.
\end{enumerate}
 \end{prop}
 
  \begin{pf}  Let $d=\deg P$; then $\syz(F, G,H)$ and $\syz(PF,PG,PH)(d)$ coincide,  and that gives the first identity. For the second identity, note that there is an injective map $\syz(F,G,H)\to \syz(PF,PG,H)(d)$ that sends  
$(\alpha,\beta, \gamma)$ to $ (\alpha,\beta, P\gamma)$. If  $P$ is prime to  $H$   then this map  is surjective as well; so   $m(PF,PG,H)= m(F,G,H)+d$, and     the identity follows easily.
\qed\end{pf}

\begin{prop}\label{lipschitz}
If $P\in \Bbbk[x,y]$ is a nonzero  homogeneous polynomial, then 
\[\vert\delta(F,G,PH)-\delta(F,G,H)\vert\le \deg P.\]
\end{prop}

\begin{pf}
Let $d=\deg P$. There is a map  $ \syz(F,G,H)\to \syz(F,G,PH)(d)$,  
$(\alpha,\beta,\gamma)\mapsto (\alpha P,\beta P,\gamma )$; so $m(F,G,PH)\le m(F,G,H)+d$.  There is also a degree-preserving map  $\syz(F,G,PH)\to \syz(F,G,H)$,  $(\alpha,\beta,\gamma)\mapsto (\alpha,\beta,\gamma P)$; so  $m(F,G,H)\le m(F,G,PH)$.   The desired inequality   follows  at once.
\qed
\end{pf}

If  $\ell\in \Bbbk[x,y]$ is a linear form, $\delta(F,G,H)$ and $\delta(F,G,H\ell)$ cannot be equal, since they have different parities. So, by the previous proposition, 
\[\delta(F,G,\ell H)=\delta(F,G,H)\pm 1.\]
We can make this   more precise: 

\begin{prop} \label{leminha}
 Suppose $F$, $G$, and $H$ have no common factor, and let  $\ell\in \Bbbk[x,y]$ be a linear form.   If $\delta(F,G,H)>0$ and $(\alpha,\beta,\gamma)$ is a syzygy of $(F,G,H)$ of minimal degree, then $\delta(F,G,\ell H)=\delta(F,G,H)+ 1$ if    $\ell$ divides $\gamma$, and $\delta(F,G,\ell H)=\delta(F,G,H)- 1$ otherwise. 
\end{prop}

\begin{pf}
If $\ell$ divides $\gamma$, then $(\alpha,\beta,\gamma/\ell)$ is an element of $\syz(F,G,\ell H)$ of minimal degree; so $m(F,G,\ell H)=m(F,G,H)$,   giving  $\delta(F,G,\ell H)=\delta(F,G,H)+1$. If  $\ell$  does not divide  $\gamma$, then we claim that $(\alpha\ell,\beta\ell,\gamma)$ is an element of  $ \syz(F,G,\ell H)$ of minimal degree. In fact, suppose there  exists $(\alpha',\beta',\gamma')\in \syz(F,G,\ell H)$ of degree $m=m(F,G,H)$. Then $(\alpha',\beta',\gamma'\ell)\in \syz(F,G, H)$ has degree $m$, and since the syzygy gap  of $F$, $G$, and $H$ is nonzero, $(\alpha',\beta',\gamma'\ell)$ must be a constant multiple of $(\alpha,\beta,\gamma)$, contradicting the assumption  that $\ell$ does not divide $\gamma$. So $m(F,G,\ell H)=m(F,G,H)+1$, and $\delta(F,G,\ell H)=\delta(F,G,H)-1$.
\qed
\end{pf}

 \begin{prop}\label{linconv}
  Let $\ell\in \Bbbk[x,y]$ be a linear form, and 
suppose $F$, $G$, and $\ell H$ have no common factor. 
If $\delta(F,G,H)$ and $\delta(F,G,\ell^{2}H)$ are both greater than $\delta(F,G,\ell 
H)$, then $\delta(F,G,\ell H)=0$.  
\end{prop}

 \begin{pf}
Multiplication by $\ell$ gives us a surjective map
\[\ideal{F,G,H}/\ideal{F,G,\ell H}\ \stackrel{\ell}{\longrightarrow}\ 
\ideal{F,G,\ell H}/\ideal{F,G,\ell^{2} H},\]
so  $\deg\ideal{F,G,\ell H}-\deg\ideal{F,G,H}\ge \deg\ideal{F,G,\ell^{2} 
H}-\deg\ideal{F,G,\ell H}$. Using Proposition \ref{formula1}   we obtain 
\[\delta(F,G,H)^{2}+\delta(F,G,\ell^{2} H)^{2}\le 2\cdot \delta(F,G,\ell H)^{2}+2.\]
But $\delta(F,G,H)=\delta(F,G,\ell^{2} H)=\delta(F,G,\ell H)+1$, so the inequality   above implies 
that $\delta(F,G,\ell H)=0$.  
\qed
\end{pf}

\begin{prop}\label{2Dconv}
Let $\ell_{1}$ 
and $\ell_{2}$ be relatively prime linear forms, such that $F$, $G$ 
and $H\ell_{1}\ell_{2}$ have no common factor. Suppose that 
$\delta(F,G,H)=\delta(F,G,H\ell_{1}\ell_{2})$ and $\delta(F,G,H\ell_{1})=\delta(F,G,H\ell_{2})$. 
Then either $\delta(F,G,H)=0$ or $\delta(F,G,H\ell_{1})=0$.
\end{prop}

\begin{pf}
Suppose $\delta:=\delta(F,G,H)>0$, 
and let  $(\alpha, \beta,\gamma)$ be a syzygy 
of  $(F, G, H)$ of minimal degree $m$. We use Proposition  \ref{leminha} repeatedly. Since
$\delta(F,G,H\ell_{1})=\delta(F,G,H\ell_{2})$, either both $\ell_1$ and $\ell_2$ divide $\gamma$, or neither one does. If both linear forms divided $\gamma$, then   $(\alpha,\beta,\gamma/(\ell_{1}\ell_{2}))$ would be a 
syzygy of  $(F,G, H\ell_{1}\ell_{2})$ of degree $m$ and we would have $\delta(F,G,H\ell_1\ell_2)>\delta$, contradicting our hypothesis.  So  neither $\ell_1$ nor $\ell_2$ divides  $\gamma$, and $\delta(F,G,H\ell_1)=\delta(F,G,H\ell_2)=\delta-1.$

Now $(\ell_{1}\alpha,\ell_{1}\beta,\gamma)$ is a syzygy of  $(F ,G, H\ell_{1})$ of minimal degree,  and since $\ell_{2}$ does not divide $\gamma$  it must be the case that  $\delta(F,G,H\ell_{1})=0$, since otherwise   $\delta(F,G,H\ell_1\ell_2)=\delta-2$, contradicting the hypothesis.
\qed\end{pf}

\section{Syzygy gap fractals}\label{gaps-and-fractals}

The  properties of syzygy gaps so far discussed  hold over  arbitrary fields. In this section, and in the remainder of the paper, we assume that $\charac \Bbbk=p>0$ and introduce a family of functions defined in terms of syzygy gaps.  Once again $F$, $G$, and $H$ are nonzero homogeneous polynomials in $A=\Bbbk[x,y]$.
If $(\alpha,\beta,\gamma)$ is a syzygy of $(F,G,H)$ of minimal degree, then $(\alpha^p,\beta^p,\gamma^p)$ is a syzygy of $(F^p,G^p,H^p)$ of minimal degree. It follows that
\begin{equation}\label{delta of p-th powers}
\delta(F^p,G^p,H^p)=p\cdot \delta(F,G,H). 
\end{equation}

In what follows, we fix a positive integer $n$ and pairwise prime linear forms  $\ell_1,\ldots,\ell_n\in \Bbbk[x,y]$. For ease of notation we introduce the following shorthands, which will be used throughout the paper:   $\ell=\prod_{i=1}^n\ell_i$, and  for any nonnegative integer vector $\vv{a}=(a_1,\ldots, a_n)$, $\ell^\vv{a}=\prod_{i=1}^n\ell_i^{a_i}$.

\begin{defn} A \emph{cell}  (with respect to the linear forms $\ell_1,\ldots,\ell_n$) is a triple $\cell{F,G,H}$ of nonzero homogeneous polynomials in $\Bbbk[x,y]$ such that $F$, $G$, and $H\ell$ have no common factor.
\end{defn}

Let $C=\cell{F,G,H}$ be a cell, $[q]=\{0,1,\ldots,q\}$, and $\vv{a}\in [q]^n$; we wish to understand how $\delta(F^q,G^q,H^q\ell^\vv{a})$ depends on $q$ and $\vv{a}$.
Eq. (\ref{delta of p-th powers}) allows us to conveniently encode these syzygy gaps   in a single function $\mathscr{I}^n\to \mathbb{Q}$, where $\mathscr{I}=[0,1]\cap\mathbb{Z}[1/p]$:
 
\begin{defn}  To  each  cell  $C=\cell{F,G,H}$ we attach a function 
$\delta_C:\mathscr{I}^n\to \mathbb{Q}$ where 
\[\delta_C\left(\frac{\vv{a}}{q}\right)=\frac{1}{q}\cdot \delta\left(F^q,G^q,H^q\ell^\vv{a}\right)\]
for any $q$ and any  $\vv{a}\in [q]^n$. (Eq. (\ref{delta of p-th powers}) ensures that  $\delta_C$ is  well-defined.) We shall nickname these functions \emph{syzygy gap fractals}, for reasons that will soon become apparent. 
\end{defn}

 \begin{rem}
In \cite{pfractals1} Monsky and the author studied a closely related  family of functions $\varphi_I$ associated to   zero-dimensional ideals $I$ of $\Bbbk \llbracket x,y\rrbracket$.  In what follows we shall explore this relation. 
\end{rem}

Let  $C=\cell{F,G,H}$ be a   cell and  $I=\ideal{F,G}:H$.  Since $\ideal{F,G,\ell}\subseteq \ideal{I,\ell}$ and $F$, $G$, and $\ell$ have no common factor, $\deg\ideal{I,\ell}<\infty$ and we can define the following function, as in \cite{pfractals1}:
\begin{align*}
\varphi_I:\mathscr{I}^n&\longrightarrow \mathbb{Q} \\
\frac{\vv{a}}{q}&\longmapsto  \frac{1}{q^2}\cdot\deg\ideal{I^{[q]},\ell^\vv{a}}
\end{align*}
Here  $I^{[q]}$ denotes  the $q$th Frobenius power of $I$, \ie the ideal generated by the $q$th powers of the elements of $I$. To relate $\delta_C$ and $\varphi_I$  
we define  a similar function 
\begin{align*}
\varphi_C:\mathscr{I}^n&\longrightarrow \mathbb{Q} \\
\frac{\vv{a}}{q}&\longmapsto \frac{1}{q^2}\cdot\deg\ideal{F^q,G^q,H^q\ell^\vv{a}} 
\end{align*} 
and start by relating $\varphi_C$ and $\delta_C$. 
Setting $d=\deg\ideal{F,G,H}$ and $\delta=\delta(F,G,H)$, Proposition \ref{formula1} gives
\begin{equation}\label{phi_C and delta_C}
 4\varphi_C(\vv{t})=\delta_C^2(\vv{t})+4d-\delta^2+2(\deg F+\deg G-\deg H)\sum_{i=1}^n{t_i}-\left(\sum_{i=1}^n t_i\right)^2.
\end{equation}
To relate $\varphi_I$ and $\varphi_C$, note that  for any ideal  $J$  of $A$ and $f\in A$  we have $\deg J=\deg(J:f)+\deg\ideal{J,f}$, and  
replacing $J$ with $\ideal{J,fg}$, that  becomes
$\deg\ideal{J,fg}= \deg\ideal{(J:f),g}+\deg\ideal{J,f}.$ Setting $J=\ideal{F^q,G^q}$, $f=H^q$, $g=\ell^\vv{a}$, and dividing by $q^2$ we find that $\varphi_C=\varphi_I+\deg\ideal{F,G,H}$. Together with Eq. (\ref{phi_C and delta_C}), this gives
\begin{equation}\label{phi_I and delta_C}
4\varphi_I(\vv{t})=\delta_C^2(\vv{t})-\delta^2+2(\deg F+\deg G-\deg H)\sum_{i=1}^n{t_i}-\left(\sum_{i=1}^n t_i\right)^2.
\end{equation}
 This, in turn, gives us the following result:

\begin{prop}\label{colon reduction} Let $C=\cell{F,G,H}$ be a   cell. Then the  ideal $\ideal{F,G}:H$ is generated by two homogeneous polynomials $U$ and $V$ such that  $\delta_C=\delta_\cell{U,V,1}$.    
\end{prop}

\begin{pf}   $\syz(F,G,H)$ has two homogeneous generators, and their third components $U$ and $V$ generate $\ideal{F,G}:H$.  Since $\ideal{F,G}\subseteq \ideal{U,V}$, the polynomials $U$, $V$, and $\ell$ have no common factor, so $\cell{U,V,1}$ is a cell.    The    generators of $\syz(F,G,H)$ have degrees $\deg U+\deg H$ and $\deg V+\deg H$, so Eq. (\ref{relation degrees vs pol degrees}) shows that $\deg U+\deg V=\deg F+\deg G-\deg H$. Noting that $\delta(U,V,1)= |\deg U-\deg V|=\delta(F,G,H)$, the result  is obtained by  replacing $C=\cell{F,G,H}$ with $\cell{U,V,1}$ in Eq.  (\ref{phi_I and delta_C}) and comparing with the  same equation in its original form.
\qed
\end{pf}

\begin{rem}
If the image of the colon ideal $I=\ideal{F,G}:H$ in $\bar{A}=A/\ideal{\ell}$ is not principal, then $\delta_C=\delta_\cell{U,V,1}$ for \emph{any} pair of generators $U$ and $V$ of $I$. This is not the case otherwise. In fact, if the image of $\ideal{U,V}$ is principal in $\bar{A}$, suppose the image of $U$ is the generator. We can modify $V$ by a multiple of $U$, without affecting $\delta_\cell{U,V,1}$, to assume that     $V=W\ell$ for some $W$.  Then Proposition \ref{deltaproperties1} shows that $\delta_\cell{U,V,1}(\vv{a}/q)=
q^{-1} \cdot \delta(U^{q},W^{q}\ell^q,\ell^\vv{a})=q^{-1}\cdot \delta(U^{q},W^{q}\ell^q/\ell^\vv{a},1)=|\deg V-\deg U-\sum_{i=1}^n a_{i}/q|$, 
which depends on the degree of $V$.   
\end{rem}

 \begin{rem}
 In view of Proposition \ref{colon reduction},  as far as the study of the functions $\delta_C$ is concerned we can always assume that the cells have the form $\cell{F,G,1}$, which we shall often abbreviate by $\cell{F,G}$. 
\end{rem}
 
 \begin{exmp}\label{r=2}
We use the above remark to explicitly describe the $\delta_C$ when  $n=2$. Suppose $C=\cell{F,G}$ is a   cell, where  $\deg F\le \deg G$. A change of variables allows us to assume that $\ell_1=x$ and $\ell_2=y$. Several cases must be considered, depending on whether or not each of $x$ and $y$ divides each of $F$ and  $G$. Suppose for instance that  $x$ divides $F$, but $y$ does not. Modifying $G$ by a multiple  of $F$, if necessary,  we can assume that $y$ divides $G$, and  for any   $a,b\le q$ we have
 \[\delta(F^q,G^q,x^{a}y^{b})=\delta(F^q/x^a,G^q/y^{b},1)=|q\deg G-q\deg F+a-b|,
\]
 by Proposition \ref{deltaproperties1}. 
Dividing by $q$ and noting that $\deg G-\deg F=\delta_C(\vv{0})$ we find 
\[
\delta_C(t_1,t_2)=|\delta_C(\vv{0})+t_1-t_2|.\]

 In all other cases similar calculations show that $\delta_C$ is a piecewise linear function of the form
 \[
 \delta_C(t_1,t_2)=|\delta_C(\vv{0})\pm t_1\pm t_2|.
 \]
The case $n=1$ is, of course, just as simple---setting $t_2=0$ in the above formula we see that $\delta_C$ is of the form
\[\delta_C(t)=|\delta_C(0)\pm t|.\]
\end{exmp}

\begin{exmp} 
In contrast, the case $n=3$ already shows some surprises. Consider for example the function $\delta_\cell{x,y}$. A linear change of variables allows us to assume that $\ell_1=x$, $\ell_2=y$, and $\ell_3=x+y$, while fixing the ideal $\ideal{x,y}$. Because of Proposition \ref{deltaproperties1}, $\delta(x^q,y^q,x^{a_1}y^{a_2}(x+y)^{a_3})=\delta(x^{q-a_1},y^{q-a_2},(x+y)^{a_3})$, so we might as well study the function 
\[
\frac{\vv{a}}{q}\longmapsto \frac{1}{q}\cdot \delta(x^{a_1},y^{a_2},(x+y)^{a_3}),
\]  
a ``reflection'' of $\delta_\cell{x,y}$. This function was studied and completely described by Han in her thesis \cite{han thesis}.  It is a Lipschitz function---a consequence of Proposition \ref{lipschitz}---and therefore can be extended (uniquely)  to a continuous function $\delta^*: [0,1]^3\to \mathbb{R}$. If $t_i>t_j+t_k$, where $\{i,j,k\}=\{1,2,3\}$,  Remark \ref{triangle-inequalities} shows that $\delta^*(\vv{t})=t_i-t_j-t_k$. If, on the other hand, the coordinates of $\vv{t}$ satisfy the triangle inequalities $t_i\le  t_j+t_k$, the description of $\delta^*(\vv{t})$ is more subtle.  Let $L_\text{odd}$ denote the elements of $\mathbb{Z}^3$ whose coordinate sum is odd, and let $d:\mathbb{R}^3\times \mathbb{R}^3 \to \mathbb{R}$ be the ``taxi-cab'' metric, $d(\vv{u},\vv{v})=\sum_{i=1}^3
|u_i-v_i|$. Then $\delta^*$ can be described as follows: 
\begin{thm}[Han \cite{han thesis}]  Suppose the coordinates of $\vv{t}\in [0,1]^3$ satisfy the triangle inequalities. If there is a pair $(s,\vv{u})\in \mathbb{Z}\times L_\text{odd}$ such that $d(p^s\vv{t},\vv{u})<1$, then there is a unique such pair with $s$ minimal. For this pair we have \[\delta^*(\vv{t})=p^{-s}(1-d(p^s\vv{t},\vv{u})).\]  If no pair $(s,\vv{u})$  
exists, then $\delta^*(\vv{t})=0$.
\end{thm}
A proof of the above result can also be found in \cite[Corollary 23]{mason}.  Figure \ref{han fractal} shows   the  two-dimensional ``slice'' $(t_1,t_2)\mapsto \delta^*(t_1,t_2,t_2)$,  where $\charac \Bbbk=3$, in the form of a relief plot,  where the color encodes the value of the function at 
each point---the higher the value, the lighter the color.
\begin{figure}
\centering
\includegraphics[height=3.5in]{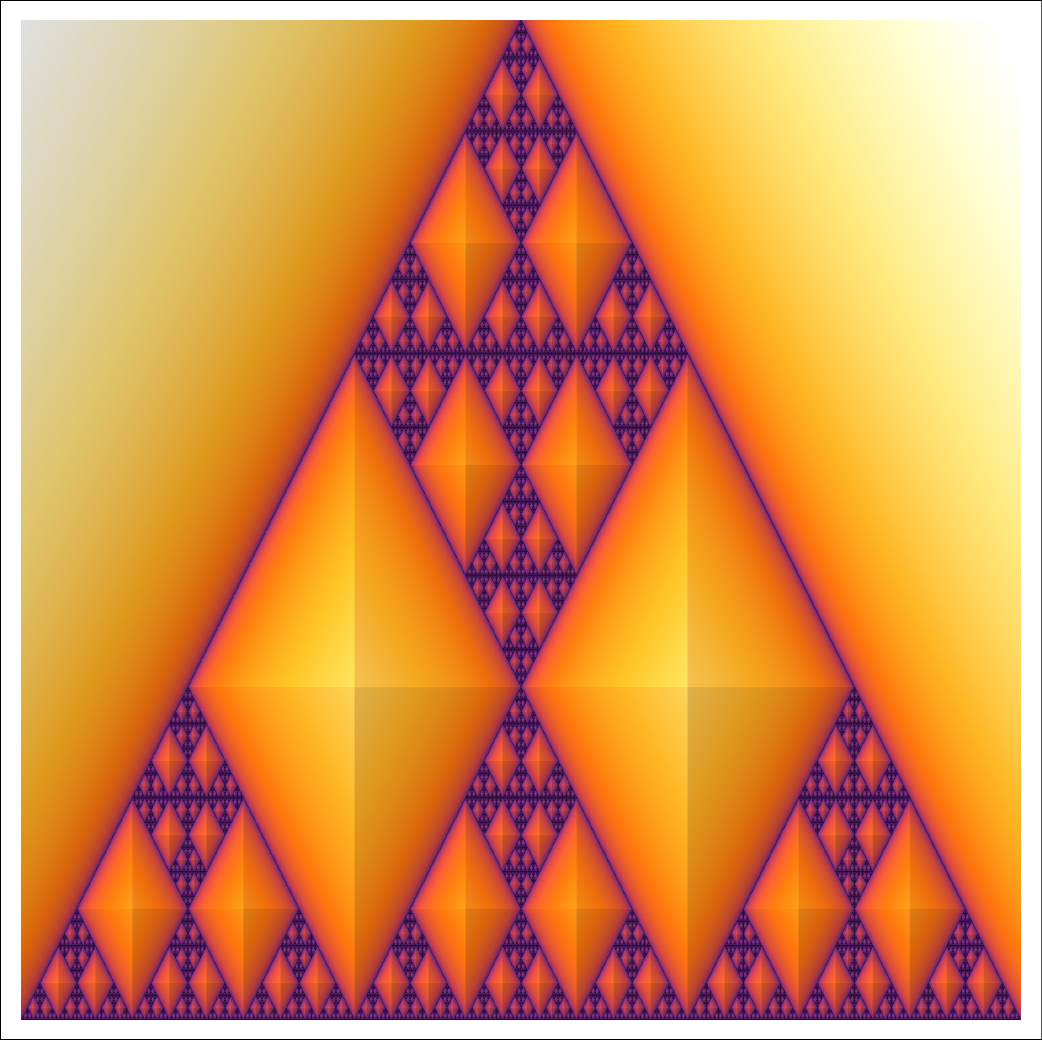} 
\caption{A two-dimensional ``slice'' of Han's  fractal $\delta^*$ in characteristic 3}
\label{han fractal}
\end{figure}
\end{exmp}

 We turn now to a couple of (related)   examples with $n=4$. 
 
 \begin{exmp}
 \label{example1D} 
Let $\Bbbk={\mathbb F}_9$, and  $\epsilon\in \Bbbk$  with $\epsilon^2+2\epsilon+2=0$; let $\ell_1,\ldots,\ell_4$ be $x$, $y$, $x+y$, and $x+\epsilon y$,  and $C=\cell{x,y}$.  We examine     the restriction of $\delta_C$ to the diagonal, namely the map  $\diag:\mathscr{I}\to \mathbb{Q}$, $\diag(t)=\delta_C(t,t,t,t)$.  
The graph of  $\diag$  is shown in Figure \ref{graph1D}. 
\begin{figure}
\centering
\includegraphics[width=3.5in]{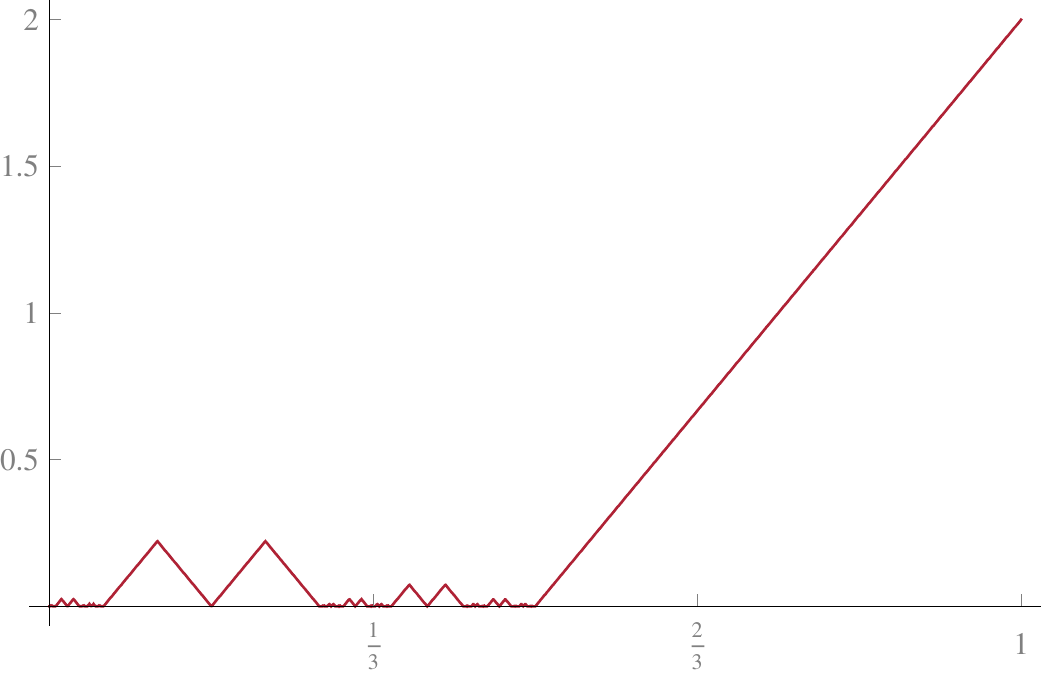}
\caption{
$\diag(t)$ $(0\le t \le 1)$}
\label{graph1D}
\end{figure}

The linear behavior on
$[1/2,1]$  is expected from Remark \ref{triangle-inequalities}: if $a/q\ge 1/2$ then $\deg(x^ay^a(x+y)^a(x+\epsilon y)^a)\ge 
2q$, so
\[\diag\left(\frac{a}{q}\right)=\frac{1}{q}\cdot \delta(x^{q},y^{q},x^ay^a(x+y)^a(x+\epsilon y)^a)=\frac{1}{q}(4a-2q)=
4\cdot \frac{a}{q}-2.\]
Note how the portion of the graph on the interval   $[1/3,2/3]$ seems to be a miniature of the entire graph. A  closer look at the portion over  $[0,1/3]$ (Figure \ref{graph1Dzoom})
\begin{figure}
\centering
\includegraphics[width=3.5in]{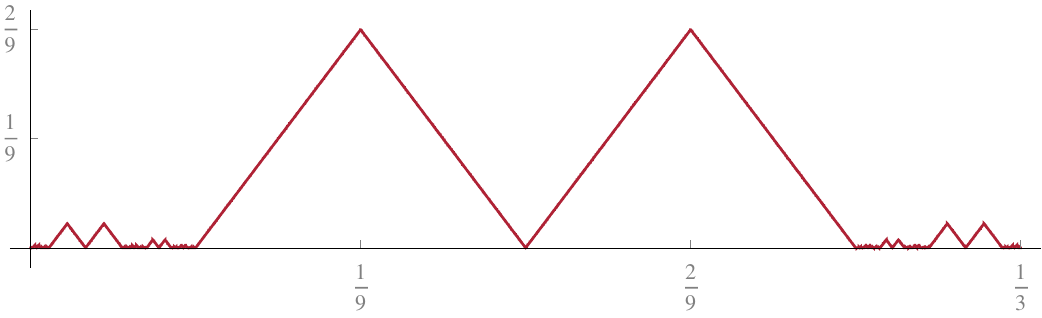} 
\caption{ 
$\diag(t)$ $(0\le t \le 1/3)$}
\label{graph1Dzoom}
\end{figure}
shows small copies of the graph of $\diag$ and its reflection about a vertical axis. These self-similarity properties will be investigated closely in a sequel  to this  paper.  

The following property can also be inferred from the graphs: at any   $t$, $\diag(t)$ seems to be simply 4 times the distance from $t$ to  the nearest zero of $\diag$---so apparently  $\diag$ can be completely reconstructed from its zeros. This is in fact the case; see Section \ref{zeros}.
\end{exmp}

\begin{exmp}\label{example2D}
With $\Bbbk$, $\ell_1,\ldots,\ell_4$, and $C$ as in the previous example, we now examine a two-dimensional ``slice'' of $\delta_C$, namely the map  $\slice:\mathscr{I}^2\to \mathbb{Q}$, $\slice(t_1,t_2)=\delta_C(t_1,t_1,t_2,t_2)$. 
 A relief plot of $\slice$ is shown in Figure \ref{delta}. 
   \begin{figure} 
\centering
\includegraphics[height=3.5in]{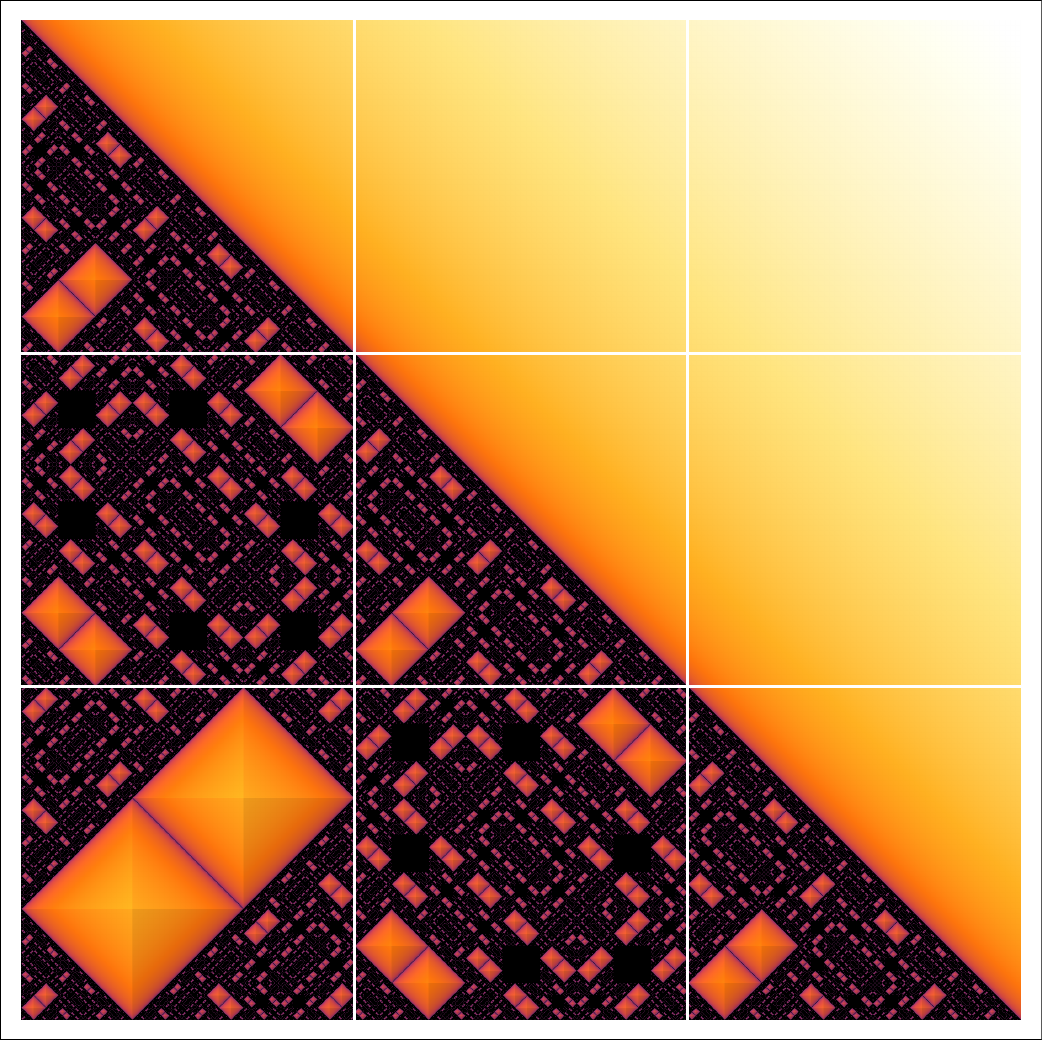}  
\caption{A two-dimensional ``slice'' of a 4-variable syzygy gap fractal}\label{delta}
\end{figure}
We immediately observe    a simple behavior on a large portion of the domain: $\slice$ is linear for $t_1+t_2\ge 1$, as expected from Remark \ref{triangle-inequalities}.  The grid dividing the plot into nine squares of equal size makes some    self-similarity properties of $\slice$ quite  evident.  (Figure \ref{fig: alpha} shows a magnification of one of those pieces---a two-dimensional ``slice'' of $\delta_\cell{x^3,y^3,xy}$, as will become clear after Section \ref{blow-up}.)
While in Section \ref{blow-up} we  discuss a couple of these self-similarity properties, their thorough study will be left for a sequel to this paper,  where we  shall  develop  the tools to verify them rigorously. 

Figure \ref{delta discrete} shows some numerical values of the function 
\[(i,j)\mapsto \frac{1}{2}\cdot\delta(x^{81},y^{81},x^iy^i(x+y)^j(x+\epsilon y)^j)=\frac{81}{2}\cdot \slice\left(\frac{i}{81},\frac{j}{81}\right),\] where zeros are replaced by dots and  the linear portion of the function is omitted.  
 \begin{figure} 
\centering
\includegraphics[height=4.77in]{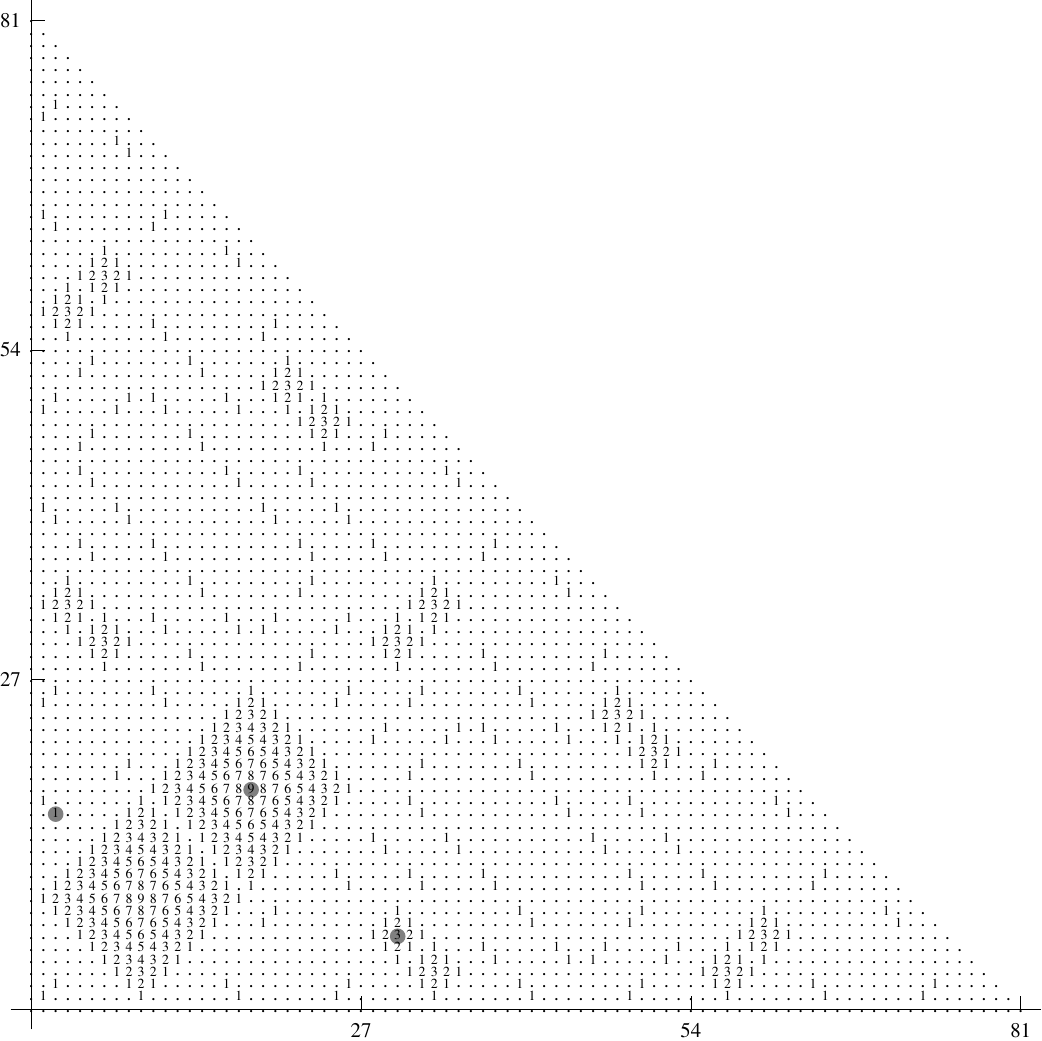} 
\caption{The function $(i,j)\mapsto  \frac{1}{2}\cdot\delta(x^{81},y^{81},x^iy^i(x+y)^j(x+\epsilon y)^j)$}\label{delta discrete}
\end{figure}
From these numerical values we infer  a property similar to that  noticed   in Example \ref{example1D}: at any point $\vv{t}$, $\slice(\vv{t})$ is simply twice the taxi-cab distance   from $\vv{t}$ to the nearest zero of $\slice$. 
\end{exmp}

\section{Syzygy gap fractals are determined by their zeros} \label{zeros}

Throughout this section we fix a cell $C$ with respect  to pairwise prime linear forms $\ell_1,\ldots,\ell_n$; 
we   use   results  from  Section \ref{gaps-and-delta}   to  prove 
 that  $\delta_C$, if nonlinear,  is completely determined by  its  zeros, as suggested  by the   examples in the previous section.  An important  role is played by  the Lipschitz property for  $\delta_C$, which follows directly from Proposition \ref{lipschitz}:

\begin{prop}\label{lip prop}
For each  $\vv{t}$ and $\vv{u}$ in $\mathscr{I}^{n}$ we have  \[|\delta_C(\vv{t})-\delta_C(\vv{u})|\le d(\vv{t},\vv{u}),\] 
where $d(\vv{t},\vv{u})$ is the \emph{taxi cab} distance between $\vv{t}$ and $\vv{u}$, 
$d(\vv{t},\vv{u})=\sum_{i=1}^{n}|t_{i}-u_{i}|.$
\end{prop}

 \begin{rem}\label{extension by continuity} 
A consequence of this is  that we can extend (uniquely)  $\delta_C$ to  a continuous function $\delta^*_C:[0,1]^n\to \mathbb{R}$. The results of this section  apply     to  $\delta^*_C$  as well, by continuity.  This extension  will be necessary in Section \ref{bounds}.
\end{rem}

In what follows,  $\Xset_{q}$ is the subset of $\mathscr{I}^{n}$ consisting of points that can be written as $(a_1/q,\ldots,a_n/q)$, with 
  $a_i\in \mathbb{Z}$. In particular, 
  $\Xset_{1}$ is simply $\{0,1\}^{n}$, the set of 
  \emph{corners} of $\mathscr{I}^{n}$.

\begin{lem}\label{distform} Suppose $\delta_C|_{\Xset_{1}}$ attains 
a ``local minimum'' at a corner $\vv{u}$, in the sense that the values 
of $\delta_C$ at all corners adjacent to $\vv{u}$ are greater than $\delta_C(\vv{u})$.
Moreover, suppose $\delta_C(\vv{u})>0$.  Then 
\begin{equation}\label{star}
 \delta_C(\vv{t})=\delta_C(\vv{u})+d(\vv{t},\vv{u}),
 \end{equation}
for all $\vv{t}\in \mathscr{I}^{n}$. In particular $\delta_C$ is linear, everywhere positive, and has a minimum at $\vv{u}$ in the usual sense.
\end{lem}

\begin{pf}
In view of our local minimum assumption,  Proposition  \ref{leminha} shows that $\delta_C(\vv{v})=\delta_C(\vv{u})+1$ for any 
corner $\vv{v}$ adjacent to $\vv{u}$. It follows from Proposition  \ref{lip prop}  that 
(\ref{star}) holds for all  points  $\vv{t}$ along  the edges containing  $\vv{u}$.

 Aiming at a contradiction, suppose that  (\ref{star})  fails for one or more points of $\Xset_{q}$. Among all such 
points, choose  $\vv{t}$ whose distance to $\vv{u}$ is minimal.  We know that $\vv{t}$ does not lie in any of the edges 
connecting to $\vv{u}$, so at least two coordinates of $\vv{t}$ and $\vv{u}$ must be different; 
say $u_{i}\ne t_{i}$ and $u_{j}\ne t_{j}$. Altering the $i$th or $j$th coordinates of $\vv{t}$ by $1/q$   we  obtain points 
$\vv{v}$, $\vv{w}$, and $\vv{z}$ that are closer to $\vv{u}$, as illustrated in Figure \ref{fig: diagram.tex}.
\begin{figure}
\centering
\unitlength .8mm
\begin{picture}(47.50,35.00)(0,0)

\put(33.00,22.20){\makebox(0,0)[cc]{$\vv{t}$}}
\put(18.00,22.00){\makebox(0,0)[cc]{$\vv{v}$}}
\put(33.13,8.00){\makebox(0,0)[cc]{$\vv{w}$}}
\put(17.50,8.00){\makebox(0,0)[cc]{$\vv{z}$}}
\put(-3.00,-2.00){\makebox(0,0)[cc]{$\vv{u}$}}

\linethickness{0.15mm}

\put(0.00,0.00){\circle*{1.75}}
\put(20.00,10.00){\circle*{1.25}}
\put(30.00,10.00){\circle*{1.25}}
\put(20.00,20.00){\circle*{1.25}}
\put(30.00,20.00){\circle*{1.25}}

\multiput(0.00,0.00)(1.95,0){21}{\line(1,0){0.98}}
\multiput(0.00,0.00)(0,1.94){16}{\line(0,1){0.97}}
\multiput(40.00,0.00)(0,1.94){16}{\line(0,1){0.97}}
\multiput(0.00,30.00)(1.95,0){21}{\line(1,0){0.98}}
\multiput(10.00,0.00)(0,1.94){16}{\line(0,1){0.97}}
\multiput(20.00,0.00)(0,1.94){16}{\line(0,1){0.97}}
\multiput(30.00,0.00)(0,1.94){16}{\line(0,1){0.97}}
\multiput(0.00,20.00)(1.95,0){21}{\line(1,0){0.98}}
\multiput(0.00,10.00)(1.95,0){21}{\line(1,0){0.98}}
\linethickness{0.40mm}
\put(0.00,0.00){\line(1,0){45.00}}
\put(0.00,0.00){\line(0,1){35.00}}

\put(32,15){\makebox(0,0)[cc]{\small{$\frac{1}{q}$}}}
\put(25,23.5){\makebox(0,0)[cc]{\small{$\frac{1}{q}$}}}
 
\put(37,-3){\makebox(0,0)[cc]{\small{$i$th coord}}}

\put(-3,26.5){\makebox(0,0)[cc]{\begin{sideways}\small{$j$th coord}\end{sideways}}}

\end{picture}

\medskip
\caption{}\label{fig: diagram.tex}
\end{figure}
Let $s=\delta_C(\vv{v})$; because of our choice of $\vv{t}$,  we can use  (\ref{star})    to conclude that 
$\delta_C(\vv{w})=s$ 
and $\delta_C(\vv{z})=s-\frac{1}{q}$. But, since $\delta_C(\vv{z}) \ge \delta_C (\vv{u}) >0$, Proposition  \ref{2Dconv} then says  that $\delta_C(\vv{t})$ must equal $s+\frac{1}{q}$, and (\ref{star}) holds for $\vv{t}$, 
contradicting our assumption. 
\qed
\end{pf}

\begin{bigthm}  \label{thm: delta det by 0s}
 Let $\ZEROset=\{\vv{z}\in \mathscr{I}^{n}\,|\, \delta_C(\vv{z})=0\}$. 
\begin{itemize}
\item If $\ZEROset$ is empty, then we are in the situation of  
Lemma \ref{distform}, and $\delta_C$ is linear; it takes on a minimum value at a corner  $\vv{u}$ of $\mathscr{I}^n$ and,  at each     $\vv{t}\in \mathscr{I}^n$, \[\delta_C(\vv{t})=\delta_C(\vv{u})+d(\vv{t},\vv{u}).\]
\item
If $\ZEROset$ is nonempty then $\delta_C(\vv{t})$ is the  taxi-cab distance from
$\vv{t}$ to the set $\ZEROset$, for all $\vv{t}\in \mathscr{I}^{n}$.
\end{itemize}
\end{bigthm}

\begin{pf}
If $\ZEROset$ is empty, then there must be a corner $\vv{u}$ satisfying the 
hypothesis of the previous lemma. 
Suppose  $\ZEROset$ is nonempty. We shall  show the result for $\vv{t}\in \Xset_{q}$, by induction on $\psi(\vv{t})=q\cdot \delta_C(\vv{t})$. If 
$\psi(\vv{t})=0$, then $\vv{t}\in \ZEROset$  and there is nothing to show, so suppose $\psi(\vv{t})>0$.  
We claim that there is a $\vv{u}\in \Xset_{q}$ adjacent to $\vv{t}$ with 
$\psi(\vv{u})<\psi(\vv{t})$, so induction gives us the desired result. 

To prove this claim,  aiming at a contradiction suppose that $\delta_C(\vv{u})>\delta_C(\vv{t})>0$, 
for all $\vv{u}$ adjacent to $\vv{t}$ in $\Xset_q$. Then Proposition \ref{linconv} shows 
that $\vv{t}$ must be a corner. If there were some adjacent corner $\vv{v}$ with $\delta_C(\vv{v})<\delta_C(\vv{t})$, then $\delta_C(\vv{t})-\delta_C(\vv{v})=1=d(\vv{t},\vv{v})$, and Proposition \ref{lip prop} would show that $\delta_C$ is linear  on the edge  linking $\vv{t}$ and $\vv{v}$ (hence decreasing as one  goes from $\vv{t}$ to $\vv{v}$); but that is not possible,  as we are assuming that  $\vv{t}$ is a local minimum in $\Xset_q$. This shows that  $\vv{t}$ satisfies the 
hypothesis of Lemma \ref{distform}.  So $\delta_C$  is  everywhere positive---but this contradicts our assumption that $\ZEROset$ is nonempty.
\qed
\end{pf}

\begin{bigthm}\label{max-corner} Suppose $\delta_C|_{\Xset_{q}}$ attains 
a ``local maximum'' at  $\vv{u}$, in the sense that the values 
of $\delta_C$ at all  points of $\Xset_q$ adjacent to $\vv{u}$ are  smaller  than $\delta_C(\vv{u})$.  Then 
\begin{equation}\label{max distance formula}
\delta_C(\vv{t})=\delta_C(\vv{u})-d(\vv{t},\vv{u}),
\end{equation}
for all $\vv{t}\in \mathscr{I}^{n}$ with $d(\vv{t},\vv{u})\le \delta_C(\vv{u})$.  In particular, $\delta_C$ is piecewise linear on that region, and has a local maximum at $\vv{u}$ in the usual sense. 
\end{bigthm}

\begin{pf} We   first  prove the theorem for points of $\Xset_q$.  If  (\ref{max distance formula})  is false  for  some $\vv{t}\in \Xset_{q}$  with $d(\vv{t},\vv{u})\le \delta_C(\vv{u})$,   choose one such   $\vv{t}$ whose distance to $\vv{u}$ is minimal.  If two or more coordinates of $\vv{t}$ and $\vv{u}$ are different, the argument used in Lemma \ref{distform} yields a contradiction; so suppose  only the $i$th coordinates of $\vv{t}$ and $\vv{u}$ differ. Proposition  \ref{leminha} and the ``local maximum'' assumption show that $\vv{t}$ cannot be adjacent to $\vv{u}$ in $\Xset_q$.  Modifying the $i$th coordinate of $\vv{t}$ by $ 1/q$ and $ 2/q$ we obtain points $\vv{v}$ and $\vv{w}$ closer to $\vv{u}$, as illustrated in Figure \ref{fig: diagram2.tex}.
\begin{figure}
\centering
\unitlength 1  mm
\begin{picture}(60,3)(0,0)
\linethickness{0.15mm}
\put(0,0){\circle*{1.5}}
\put(20,0){\circle*{1}}
\put(30,0){\circle*{1}}
\put(40,0){\circle*{1}}

\put(0,-2.6){\makebox(0,0)[cc]{$\vv{u}$}}

\put(20,-2.6){\makebox(0,0)[cc]{$\vv{w}$}}
\put(30,-2.6){\makebox(0,0)[cc]{$\vv{v}$}}
\put(40,-2.6){\makebox(0,0)[cc]{$\vv{t}$}}

\put(35,2.6){\makebox(0,0)[cc]{$1/q$}}
\put(25,2.6){\makebox(0,0)[cc]{$1/q$}}

\linethickness{0.4mm}
\put(0,0){\line(1,0){65}}
\put(58.5,-2){\makebox(0,0)[cc]{\small{$i$th coord}}}

\put(0,26){\makebox(0,0)[cc]{}}

\end{picture}
\medskip
\caption{}\label{fig: diagram2.tex}
\end{figure}
Because of  our choice of $\vv{t}$, Eq. (\ref{max distance formula}) holds for these points, so  $\delta_C(\vv{w})>\delta_C(\vv{v})>0$, and Proposition \ref{linconv} shows that $\delta_C(\vv{t})=\delta_C(\vv{v})-1/q=\delta_C(\vv{u})-d(\vv{t},\vv{u})$, contradicting our assumption.

To complete the proof, note that since (\ref{max distance formula}) holds for points adjacent to $\vv{u}$ in $\Xset_q$,   it also holds on the  edges connecting these adjacent points to $\vv{u}$, by Proposition \ref{lip prop}. This shows that $\delta_C$ restricted to any $\Xset_{q'}$ with $q'\ge q$ also has a local maximum at $\vv{u}$; thus (\ref{max distance formula}) holds for all $\vv{t}\in \Xset_{q'}$ with $d(\vv{t},\vv{u})\le \delta_C(\vv{u})$.
\qed 
\end{pf}

\section{Operators on cell classes}\label{operators}

In this section we introduce  a notion of equivalence of cells,  and   present a minimum on reflection and magnification operators on cell classes, to be used in Section \ref{bounds}. All cells here  are with respect to   an arbitrary fixed set of   pairwise prime linear forms $\ell_1,\ldots,\ell_n$, unless otherwise stated.

\subsection{Cell classes}

\begin{defn} 
Two   cells $C_1$ and $C_2$ are \emph{$\delta$-equivalent} if $\delta_{C_1}=\delta_{C_2}$.  The equivalence class of a cell $C=\cell{F,G,H}$ is denoted by $\mathscr{C}$ or $\cclass{F,G,H}$. 
\end{defn}

\begin{rem}\label{operations preserving classes} 
Several properties  follow immediately from the results from Sections \ref{gaps-and-delta}  and \ref{gaps-and-fractals}:
\begin{enumerate}
\item Proposition \ref{formula1}: an equivalent cell results from any change in  $\cell{F,G,H}$ that fixes the degrees of all the ideals $\ideal{F^q,G^q, H^q\ell^\vv{a}}$  and  the quantities   $Q(\deg F^q,\deg G^q, \deg H^q\ell^\vv{a})$. In particular:
\begin{itemize}
\item $\cclass{F,G,H}= \cclass{G,F,H}$; 
\item $\cclass{F,G,H}= \cclass{aF,bG,cH}$, for nonzero $a,b,c\in \Bbbk$;
\item $\cclass{F,G,H}=\cclass{F+UG+VH\ell,G,H}$ (and obvious variations), where $U$ and $V$ are homogeneous polynomials of appropriate degrees, provided $F+UG+VH\ell\ne 0$.
\end{itemize} 
\item Proposition \ref{colon reduction}: $\cclass{F,G,H}=\cclass{U,V}$, for some $U$ and $V$ such that $\ideal{U,V}=\ideal{F,G}:H$. 
\item Proposition \ref{deltaproperties1}: $\cclass{F,G,H}= \cclass{PF,PG,H}$ (and obvious variations),  for any   nonzero homogeneous polynomial   $P$ prime to $H\ell$.
\end{enumerate}
\end{rem}

\begin{defn}
Let  $\mathscr{C}$ be a cell class represented by a cell $C$. We define $\delta_\mathscr{C}=\delta_C$ and $\delta^*_\mathscr{C}=\delta^*_C$. 
\end{defn}

 \subsection{Reflections}
 
 Let    
$R_{1}:\mathscr{I}^n\to \mathscr{I}^n$ be the reflection in the first coordinate,  \ie the map 
  that takes $(t_{1},t_2,\ldots,t_{n})$ to 
$(1-t_{1},t_{2},\ldots,t_{n})$. Given   a cell class  $\mathscr{C}=\cclass{F,G}$, we shall   construct another cell class $R_1\mathscr{C}$ such that 
$\delta_{R_1\mathscr{C}}=\delta_{\mathscr{C}}\circ R_{1}$.

 We may choose  
 $F$ and $G$ 
 of degree $\ge n$. Modifying one of 
these elements by a multiple of the other, we may assume that 
$F=\ell_{1}F^{*}$, for some $F^*$. Modifying $F$ by a multiple of $\ell$ we may 
assume that $\ell_{1}$ does not divide $F^{*}$. Since $\deg G\ge n$, $G$ 
is congruent to some $\ell_{1}G^{*}\pmod {\ell_{2}\cdots\ell_{n}}$, with 
$G^{*}\ne 0$. It is easy to see that $F^*$, $\ell_1G^*$, and $\ell$ have no common factor. Now let $R_1\mathscr{C}=\cclass{F^{*},\ell_{1}G^{*}}.$

\begin{prop}\label{reflection-relation}
$\delta_{R_{1}\mathscr{C}}=\delta_{\mathscr{C}}\circ R_{1}$.
\end{prop}

\begin{pf}  
By Proposition \ref{deltaproperties1},
\[\delta_{\mathscr{C}}\left(\frac{\vv{a}}{q}\right) =  
\frac{1}{q}\cdot
\delta\left(\ell_{1}^{q}{F^{*}}^{q},G^{q},\prod_{i=1}^{n}\ell_{i}^{a_{i}}\right) 
=\frac{1}{q}\cdot
\delta\left(\ell_{1}^{q-a_{1}}{F^{*}}^{q},G^{q},\prod_{i=2}^{n}\ell_{i}^{a_{i}}\right).\]
  Since     $G^q-\ell_1^q {G^*}^q$ is  a multiple of  $\prod_{i=2}^n \ell_i^{a_i}$,  by design, $G^q$ may be replaced with $\ell_1^q {G^*}^q$, and a couple more uses of Proposition \ref{deltaproperties1} gives us 
\begin{align*}
\delta_{\mathscr{C}}\left(\frac{\vv{a}}{q}\right) & =  
\frac{1}{q}\cdot
\delta\left(\ell_{1}^{q-a_{1}}{F^{*}}^{q},\ell_{1}^{q}{G^{*}}^{q},\prod_{i=2}^{n}\ell_{i}^{a_{i}}\right) \\
&=\frac{1}{q}\cdot
\delta\left({F^{*}}^{q},\ell_{1}^{a_{1}}{G^{*}}^{q},\prod_{i=2}^{n}\ell_{i}^{a_{i}}\right)\\
&=\frac{1}{q}\cdot
\delta\left({F^{*}}^{q},\ell_{1}^{q}{G^{*}}^{q},\ell_{1}^{q-a_{1}}\prod_{i=2}^{n}\ell_{i}^{a_{i}}\right) \\
&=\delta_{R_1\mathscr{C}}\left(R_1\left(\frac{\vv{a}}{q}\right)\right). 
\end{align*}
\qed
\end{pf}

In particular, Proposition \ref{reflection-relation} shows that $\delta_{R_1\mathscr{C}}$ only depends on the cell class $\mathscr{C}$, and therefore the class $R_1\mathscr{C}$ is independent of the many choices made in its construction. So we have a well-defined operator $R_1$ on the set of cell classes. Furthermore, since $\delta_{R_1R_1\mathscr{C}}=\delta_{\mathscr{C}}\circ R_1\circ R_1=\delta_{\mathscr{C}}$,  it follows that $R_1R_1\mathscr{C}=\mathscr{C}$, for any class $\mathscr{C}$, so $R_1$ is an involution on the set of cell classes.  We may, of course, construct other reflection operators.

\begin{defn}
Let $\mathscr{C}$ be a cell class and $1\le i\le n$. Choose a representative $\cell{\ell_iF^*,G}$ for $\mathscr{C}$  where $\deg \ell_iF^*, \deg G\ge n$ and $\ell_i $ does not divide $ F^*$. Choose $G^*\ne 0$ such that $G\equiv \ell_iG^*\pmod{\ell/\ell_i}$. Then  we define
$R_i\mathscr{C}=\cclass{F^*,\ell_iG^*}.$
The $R_i$ are well-defined  commuting involutions on the set of cell classes, and  \[\delta_{R_i\mathscr{C}}=\delta_\mathscr{C}\circ R_i.\]  We call the $R_i$ and their compositions  \emph{reflection operators}; if $R$ is a reflection operator we call  $R\mathscr{C}$ a \emph{reflection} of $\mathscr{C}$. 
\end{defn}

For later use, we prove that  cell classes with a certain special property are unique up to reflection. 
 
\begin{defn}\label{special cell} Suppose $n\ge3$. 
 A cell class $\mathscr{C}$ is  \emph{special}
if $\delta_{\mathscr{C}}=n-2$ at some  corner $\vv{c}$ of $\mathscr{I}^{n}$, 
$\delta_{\mathscr{C}}=n-3$ at all corners adjacent to $\vv{c}$, and $\delta_{\mathscr{C}}=0$ at the  corner opposite to $\vv{c}$. 
\end{defn}

The prototypical example of a special cell class is $\cclass{x,y}$. In fact, after a change of variables we may assume that $\ell_1=x$ and $\ell_2=y$, and it is then easy to see that   $\delta_\cclass{x,y}(\vv{1})=n-2$, where $\vv{1}=(1,\ldots,1)$, and  $\delta_\cclass{x,y}=n-3$ at any corner adjacent to $\vv{1}$,  while    $\delta_\cclass{x,y}(\vv{0})=\deg x-\deg y=0$.

\begin{prop}\label{special cells are unique}
Special cell classes are reflections of one another. In particular,  each special cell class is a reflection of $\cclass{x,y}$.
 \end{prop}

\begin{pf} 
It    suffices to show that there is only one special  
cell class $\mathscr{C}$ with  $\delta_{\mathscr{C}}(\vv{0})=n-2$; any other special cell class will necessarily  be a reflection of $\mathscr{C}$. The cell class $\mathscr{C}$ may be represented by a cell $\cell{F,G}$ with   $\deg F-\deg G=\delta_{\mathscr{C}}(\vv{0})=n-2$. The assumption that  $\delta_{\mathscr{C}}=n-3$  at all corners adjacent to $\vv{0}$     implies that   
     $G$ is prime to each $\ell_i$, by Proposition \ref{leminha}. The assumption that $\delta_\mathscr{C}(\vv{1})=0$ implies that $F\not\in \ideal{G,\ell}$. (If  $F=UG+V\ell$,  with $U\ne 0$, then $\cclass{F,G}=\cclass{UG,G}=\cclass{U,1}$, and $\delta_\mathscr{C}(\vv{1})=\delta(U,1,\ell)=2$;  a similar contradiction is obtained if $U= 0$.)

Since $F\not\in \ideal{G,\ell}$,  the polynomial $G$ is not constant, so $\deg F\ge n-1$.  Thus    the $\Bbbk$-vector space of elements of degree $\deg F$ in $\Bbbk[x,y]/\ideal{\ell}$
is  $n$-dimensional.   
A basis for that space  consists of   the  elements   represented by 
\begin{equation}\label{generators}
F, \ x^{n-2}G, \ x^{n-3}yG, \ \ldots,\ y^{n-2}G,
\end{equation} which   are linearly independent   because   $F\not\in \ideal{G,\ell}$ and $G$ is prime to $\ell$.

 Now suppose $\cclass{F_{1},G_{1}}$ is 
another   cell class with the same properties.
Remark \ref{operations preserving classes}(3) allows us  to multiply  $\cell{F,G}$ and $\cell{F_1,G_1}$ by homogeneous polynomials 
prime to each $\ell_{i}$, 
so we may assume that
$G_{1}=G$ and, \emph{a fortiori}, $\deg F_1=\deg F$.  
So the image of $F_{1}$  in $\Bbbk[x,y]/\ideal{\ell}$  can be written as a linear combination of 
the images of the elements in \eqref{generators}, and the third property  in Remark \ref{operations preserving classes}(1)  ensures  that $\cclass{F_1,G}=\cclass{F,G}=\mathscr{C}$.
\qed
\end{pf}

\begin{exmp}\label{reflection of special cell}
For later use, let us find a representative for  the special cell class $\mathscr{C}=R_2\cdots R_n\cclass{x,y}$. Changing variables, if necessary, we may assume that $\ell_1=x$ and $\ell_2=y$.  Then 
\begin{align*}
\delta_\mathscr{C}\left(\frac{\vv{a}}{q}\right)&= \delta_\cclass{x,y}\left(R_2\cdots R_n\left(\frac{\vv{a}}{q}\right)\right) \\
&= \frac{1}{q}\cdot \delta\left(x^q,y^q,x^{a_1}y^{q-a_2}\ell_3^{q-a_3}\cdots \ell_n^{q-a_n}\right). 
\end{align*}
Through repeated uses of   Proposition \ref{deltaproperties1} we find that
\begin{align*}
\delta\left(x^q,y^q,x^{a_1}y^{q-a_2}\ell_3^{q-a_3}\cdots \ell_n^{q-a_n}\right)&=\delta\left(x^{q-a_1},y^{a_2},\ell_3^{q-a_3}\cdots \ell_n^{q-a_n}\right) \\
&=\delta\left(x^{q},x^{a_1}y^{a_2}\ell_3^{a_3}\cdots \ell_n^{a_n},(\ell_3\cdots\ell_n)^{q}\right), 
\end{align*}
so $\delta_\mathscr{C}(\vv{a}/q)=\delta_\cclass{x,\ell_3\cdots\ell_n}(\vv{a}/q)$, whence $\mathscr{C}=\cclass{x,\ell_3\cdots\ell_n}$. But $\ell_3\cdots\ell_n\equiv cy^{n-2}\pmod{x}$, for some nonzero constant $c$, so we conclude that 
\[R_2\cdots R_n\cclass{x,y}=\cclass{x,y^{n-2}}.\]
Of course, in view of Proposition \ref{special cells are unique}, the above identity could be just as easily verified by showing that $\cclass{x,y^{n-2}}$ is the special  cell class with maximum at $(1,0,\ldots,0)$. 
\end{exmp}

\subsection{Magnification  operators}\label{blow-up}

\begin{defn}\label{mag functions}
 Let $q$ be a power of $p$ and $\vv{b}\in [q-1]^n$. Given  $f:\mathscr{I}^n\to \mathbb{Q}$ we define $T_{q|\vv{b}}f:\mathscr{I}^n\to \mathbb{Q}$ as follows:
\[T_{q|\vv{b}}f(\vv{t})=q\cdot f\left(\frac{\vv{t}+\vv {b}}{q}\right).\]
\end{defn}

\begin{rem} This definition differs slightly from the one given in  \cite{pfractals1,pfractals2}.  
\end{rem}

We introduce operators on cell classes that are ``compatible'' with the action of the operators  $T_{q|\vv{b}}$ on the functions $\delta_\mathscr{C}$. 

\begin{defn} Let $\mathscr{C}=\cclass{F,G,H}$.
With notation as in Definition \ref{mag functions}, we define    $T_{q|\vv{b}}\mathscr{C}= \cclass{F^q,G^q,H^q\ell^{\vv{b}}}.$
Then
\[\delta_{T_{q|\vv{b}}\mathscr{C}}=T_{q|\vv{b}}\delta_\mathscr{C},\]
so $T_{q|\vv{b}}\mathscr{C}$ does not depend on the choice of representative for $\mathscr{C}$. We call $T_{q|\vv{b}}$ a \emph{magnification  operator} and $T_{q|\vv{b}}\mathscr{C}$ a \emph{magnification} of $\mathscr{C}$.
\end{defn}

\begin{exmp}\label{self-sim 1}
Suppose $\ell_1=x$ and $\ell_2=y$. Since    $\cclass{x^q,y^q,x^{j}y^{k}}=\cclass{x^{q-j},y^{q-k}}$, for $j,k\le q$, we conclude  that  
 \begin{equation}\label{blowup formula}
 T_{q|(j,k,0,\ldots,0)}\cclass{x,y}=\cclass{x^{q-j},y^{q-k}}.
 \end{equation}
 In particular, taking $q\ge n-2$ and setting $j=q-1$ and $k=q-n+2$ we find \[T_{q|(q-1,q-n+2,0,\ldots,0)}\cclass{x,y}=\cclass{x,y^{n-2}}=R_2\cdots R_n\cclass{x,y},\]
 where the second equality comes from  Example \ref{reflection of special cell}. This reveals an interesting self-similarity property of the syzygy gap fractal  $\delta_\cclass{x,y}$:
 \begin{equation}\label{eqn: self-sim}
 T_{q|(q-1,q-n+2,0,\ldots,0)}\delta_\cclass{x,y}=\delta_\cclass{x,y}\circ R_2\cdots R_n.
 \end{equation}
 
 Going back to (\ref{blowup formula}) and setting $j=k=q-1$, we see that $\cclass{x,y}$ is fixed by the operator $T_{q|(q-1,q-1,0,\ldots,0)}$ and, consequently, so is $\delta_\cclass{x,y}$.  The same holds, of course, for any  $T_{q|\vv{b}}$  where $\vv{b}$ is a permutation of $(q-1,q-1,0,\ldots,0)$.
This    self-similarity property   can be observed in Example \ref{example2D}---it explains why the NW and SE portions  of the relief plot in  Figure \ref{delta} are miniatures of the whole plot.
The fact that the center portion  is also a miniature of  the whole plot  is a consequence of the  following result (of which the  property discussed in this paragraph is a particular case, with $j=k=0$).
\end{exmp}
 
\begin{prop}Suppose $n\ge 4$ and  let $\vv{b}=(q-j-1,q-k-1,j,k,0,\ldots,0)$, for some $j,k<q$.  Let $\lambda$ be the cross ratio of the roots in $\mathbb{P}^1(\Bbbk)$ of $\ell_1$, $\ell_2$, $\ell_3$, and $\ell_4$, and suppose  $\alpha=\sum_{i=0}^k\binom{j}{k-i}\binom{k}{i}\lambda^i  \ne 0$. Then $\cclass{x,y}$  is fixed by  $T_{q|\vv{b}}$ (and consequently so is $\delta_\cclass{x,y}$). A similar result holds for any permutation of $\vv{b}$. 
\end{prop}

\begin{pf}A change of variables allows us to assume that $\ell_1=x$, $\ell_2=y$, $\ell_3=x+y$, and $\ell_4=x+\lambda y$. Then $\alpha$ is the coefficient of $x^jy^k$ in $(x+y)^j(x+\lambda y)^k$.  By Proposition \ref{deltaproperties1},  
\begin{align*}
T_{q|\vv{b}}\cclass{x,y}&=\cclass{x^q,y^q,x^{q-j-1}y^{q-k-1}(x+y)^j(x+\lambda y)^k}\\
&=\cclass{x^{j+1},y^{k+1},(x+y)^j(x+\lambda y)^k}. 
\end{align*} 
All terms of $(x+y)^j(x+\lambda y)^k$  but $\alpha x^jy^k$ are multiples of $x^{j+1}$ or $y^{k+1}$, so  $T_{q|\vv{b}}\cclass{x,y}= \cclass{x^{j+1},y^{k+1},\alpha x^jy^k}=\cclass{x,y,\alpha}=\cclass{x,y}$. 
\qed
\end{pf}

When $j=k=1$ the condition on $\alpha$ in the above proposition is simply that $\lambda\ne -1$. This is the case in Example \ref{example2D}, and explains why the center  portion of the relief plot    shown in Figure \ref{delta}  is a miniature of  whole plot. 

 We end this section with a proof that the $\delta_\mathscr{C}$ are $p$-fractals when the field $\Bbbk$ is finite. We recall the definition of $p$-fractal first:
 
\begin{defn} 
 A function  $f:\mathscr{I}^n\to \mathbb{Q}$  is a \emph{$p$-fractal} if the $\mathbb{Q}$-vector space spanned by $f$ and all the magnifications  $T_{q|\vv{b}}f$ is finite dimensional.
\end{defn}

\begin{thm}
If the field $\Bbbk$ is finite, then there are only finitely many nonlinear $\delta_\mathscr{C}$. In particular, the $\delta_\mathscr{C}$ are $p$-fractals. 
\end{thm}

\begin{pf}
We start by showing that every cell class has a representative $\cell{F,G}$ with $\deg F\le n$. In fact, 
suppose  $\mathscr{C}=\cclass{F,G}$, with  $\deg F$ and  $\deg G$ greater than $n$. Take $U$ of degree $\le 2$  and  prime   to $\ell$; it is easy to see that  $\ideal{F,G}\subseteq\ideal{x,y}^{\deg U+n-1}\subseteq \ideal{U,\ell}$, so by modifying $F$ and $G$ by multiples of $\ell$ we may assume that both are multiples of $U$. Since $U$ is prime to $\ell$ we can also divide $F$ and $G$ by $U$ without affecting $\mathscr{C}$, obtaining a new representative   consisting of polynomials of smaller degrees. 

Now note that if  $\deg G-\deg F\ge n$, Proposition \ref{deltaproperties1} and Remark \ref{triangle-inequalities}  show that $\delta_\mathscr{C}$ is linear.  Together with the result from the previous paragraph, this  shows that any cell class $\mathscr{C}$ with nonlinear $\delta_\mathscr{C}$   can be represented by a cell $\cell{F,G}$ with $\deg F\le n$ and $\deg G-\deg F< n$. If  $\Bbbk$ is finite, there are only finitely many such cells. 

The conclusion that  the $\delta_\mathscr{C}$  are $p$-fractals  follows at once,  since the $\mathbb{Q}$-vector space spanned by the finitely many nonlinear $\delta_\mathscr{C}$, the constant function 1, and the coordinate functions   is stable under the   operators $T_{q|\vv{b}}$.
\qed\end{pf}

\section{An upper bound}\label{bounds}

Throughout this section we fix   pairwise prime linear forms $\ell_1,\ldots,\ell_n$   in $\Bbbk[x,y]$. Cells and cell classes are defined with respect to these linear forms, unless otherwise stated.

In \cite[Theorem 8]{mason} Monsky  found   an upper bound  for the  local maxima of the $\delta_\mathscr{C}$: if $\delta_\mathscr{C}=\delta_\cclass{F,G,H}$ has a local maximum at $\vv{a}/q$, where  $q>1$ and some $a_i$ is not divisible by $p$, then    \[\delta_\mathscr{C}\left(\frac{\vv{a}}{q}\right)\le \frac{n_0-2}{q},\] where $n_0$ is the number of  zeros   of $FGH\ell$ in  $\mathbb{P}^1\left(\bar{\Bbbk}\right)$ (not counted with multiplicity). 

\begin{exmp} In Examples \ref{example1D} and \ref{example2D} we looked at   ``slices'' of a  syzygy gap fractal  $\delta_\mathscr{C}:\mathscr{I}^4\to  \mathbb{Q}$;  Figure \ref{delta discrete} shows some related numerical values. From  the three points highlighted   in that  picture  we can gather (with the help of Theorem \ref{max-corner}) that $\delta_\mathscr{C}$ has   
 local maxima  at the points $(2,2,2,2)/9$,  $(10,10,2,2)/27$,  and $(2,2,16,16)/81$,  where it takes on the values  $2/9$,  $2/27$, and $2/81$, respectively. Since here  $n_0=4$, Monsky's bound is attained in each case.
\end{exmp}

In this section  we sharpen Monsky's result, proving the following:

\begin{bigthm}\label{bounds1}
Suppose $\delta_\mathscr{C}$ has a local maximum at $\vv{a}/q$, where $q>1$ and $\vv{a}/q$ is reduced, in the sense  that some coordinate  $a_i/q$ is  reduced.  Then    \[\delta_\mathscr{C}\left(\frac{\vv{a}}{q}\right)\le \frac{n-2}{q}.\] 
\end{bigthm}

\begin{rem}\label{rem: bound for special cell classes} If $\mathscr{C}$ is a special cell class (see  Definition \ref{special cell}), this is nothing but Monsky's bound, since in view of Proposition \ref{special cells are unique} we may assume that  $\mathscr{C}=\cclass{x,y}=\cclass{\ell_1,\ell_2}$,  whence  $n_0=n$. 
\end{rem}

The approach used in our proof  was suggested by Monsky, and  follows closely  an alternate proof he provided of his result from \cite{mason} (private communication). Before we dive into the proof of the theorem, we look at   a couple of consequences.

\begin{cor}\label{linear-max} Let $q>1$, and fix $(a_2,\ldots,a_n)\in [q]^{n-1}$.
Suppose the map \[t\mapsto \delta_\mathscr{C}(t,a_2/q,\ldots,a_n/q)\]   has a local maximum at $a_1/q$, where $a_1$ is not divisible by $p$, and let $\vv{a}=(a_1,a_2,\ldots,a_n)$.  Then  \[\delta_\mathscr{C}\left(\frac{\vv{a}}{q}\right)\le \frac{n-2}{q}.\] 
\end{cor}

\begin{pf} According to Theorem \ref{max-corner}, each local maximum $\vv{u}$ of $\delta_\mathscr{C}|_{\Xset_{q}}$ determines a region on which 
$\delta_\mathscr{C}$ is piecewise linear, given by $\delta_\mathscr{C}(\vv{t})=\delta_\mathscr{C}(\vv{u})-d(\vv{t},\vv{u})$. The point $\vv{a}/q$ is in one such region, and because the map $t\mapsto \delta_\mathscr{C}(t,a_2/q,\ldots,a_n/q)$  has a local maximum at $a_1/q$, it must be the case that  $u_1=a_1/q$. So $\vv{u}$ satisfies the assumptions of Theorem \ref{bounds1}, and 
 $\delta_\mathscr{C}\left(\vv{a}/q\right)\le\delta_\mathscr{C}\left(\vv{u}\right)\le  (n-2)/q.$
\qed
\end{pf}

The next corollary provides an answer to  a  question   raised  in \cite[Section 7(4)]{pfractals1} in a special case.

\begin{cor} Let  $C=\cell{F,G,H}$ be a cell, and $\vv{a}=(a_1,\ldots,a_n)\in [q]^n$   with  $a_1$ not divisible by $p$.  Then 
\begin{equation}\label{newbound}2\deg\ideal{F^q,G^q,H^q\ell^\vv{a}}-\deg\ideal{F^q,G^q,H^q\ell^\vv{a}\ell_1}-\deg\ideal{F^q,G^q,H^q\ell^\vv{a}/\ell_1}\le n-2.
\end{equation}
\end{cor}

\begin{pf}   Let $\delta$, $\delta_+$, and $\delta_-$ denote the  syzygy gaps correspondent to the degrees on the  left  hand side of (\ref{newbound}), namely     $\delta(F^q,G^q,H^q\ell^\vv{a})$,
$\delta(F^q,G^q,H^q\ell^\vv{a}\ell_1)$, and 
$\delta(F^q,G^q,H^q\ell^\vv{a}/\ell_1)$.
 Proposition \ref{formula1} transforms  (\ref{newbound}) into \[\frac{1}{4}(2+2\delta^2-\delta_+^2-\delta_-^2)\le n-2.\]
If $\delta_-=\delta-1$ and $\delta_+=\delta+1$ (or vice-versa), then $2+2\delta^2-\delta_+^2-\delta_-^2=0$. If $\delta_\pm=\delta+1$, then $\delta=0$, by Proposition \ref{linconv}, so again $2+2\delta^2-\delta_+^2-\delta_-^2=0$. Finally, if $\delta_\pm=\delta-1$, then $(2+2\delta^2-\delta_+^2-\delta_-^2)/4=\delta$. But in this situation the map $t\mapsto \delta_\mathscr{C}(t,a_2/q,\ldots,a_n/q)$ has a local  maximum at $a_1/q$, and Corollary \ref{linear-max} shows that $\delta\le n-2$.
 \qed
\end{pf}

In the remainder of this section we fix   a cell class $\mathscr{C}=\cclass{U,V}$. In view of Remark \ref{operations preserving classes}(3) we may assume that $U$ and $V$ have no common factor. Since the values of our functions remain unchanged if we extend  $\Bbbk$, we may also  assume without loss of generality that $\Bbbk$ is algebraically closed.

\subsection{Some reductions, a special case,  and a proof outline}\label{reductions}

\begin{itemize}
\item (\emph{We may assume $q=p$.})
In the situation of the statement of Theorem \ref{bounds1}, write $\vv{a}=p\cdot \vv{b}+\vv{c}$, with $0\le c_i<p$, and let $q'=q/p$. Then 
\[\delta_\mathscr{C}\left(\frac{\vv{a}}{q}\right)=\delta_\mathscr{C}\left(\frac{\vv{c}/p+\vv{b}}{q'}\right)=\frac{1}{q'}\cdot (T_{q'|\vv{b}}\delta_\mathscr{C})\left(\frac{\vv{c}}{p}\right)=\frac{1}{q'}\cdot\delta_{T_{q'|\vv{b}}\mathscr{C}}\left(\frac{\vv{c}}{p}\right),\]
and  $\delta_{T_{q'|\vv{b}}\mathscr{C}}$ has a local maximum at $\vv{c}/p$.  Moreover, since $\vv{a}/q$ is reduced, so is $\vv{c}/p$. The above equation   also shows that $\delta_\mathscr{C}(\vv{a}/q)\le (n-2)/q$ whenever $\delta_{T_{q'|\vv{b}}\mathscr{C}}(\vv{c}/p)\le (n-2)/p$.   So it suffices to prove Theorem \ref{bounds1} for $q=p$. 

\item (\emph{The  trivial cases $n = 1, 2$.}) Theorem \ref{bounds1} is vacuously  true when $n=1$ or 2, since in those cases the   $\delta_{\mathscr{C}}$, completely  described  in  Example  \ref{r=2}, only have local maxima at the endpoints of $\mathscr{I}$ or corners of $\mathscr{I}^2$.

\item (\emph{Focusing on interior points.}) Note that each restriction of $\delta_\mathscr{C}$ to a face of $\mathscr{I}^n$ agrees with the values of a function $\delta_{\mathscr{C}'}:\mathscr{I}^{n-1}\to \mathbb{Q}$. Indeed, for $\epsilon=0$ or 1, $\delta_\mathscr{C}(u_1,\ldots,u_{n-1},\epsilon)=\delta_{\mathscr{C}'}(u_1,\ldots,u_{n-1})$, where $\mathscr{C'}=\cclass{U,V,\ell_n^\epsilon}$, a cell  class defined with respect to the linear forms $\ell_1,\ldots,\ell_{n-1}$. So induction on $n$ will allow us to restrict our attention to interior points of $\mathscr{I}^n$. 
\end{itemize}

These simple remarks allow us to  prove Theorem \ref{bounds1}   for $p=2$. 

\begin{potC2}
The theorem holds for $n=1$ and 2, so we let $n>2$ and argue by induction on $n$. As shown above, it suffices to consider the case $q=p=2$.  Suppose $\delta_\mathscr{C}$ has a local maximum at $\vv{t}=\vv{a}/2$, where $\vv{a}\in[2]^n$, with some  $a_i=1$. If $\vv{t}$ lies in a face of $\mathscr{I}^n$,   the observation made  above and the  induction hypothesis  give  us the desired bound. It remains to consider  $\vv{a}=(1,\ldots,1)$. Aiming at a contradiction, we suppose $\delta_\mathscr{C}(\vv{t})>(n-2)/2$.     Theorem \ref{max-corner} shows that $(1/2,\ldots,1/2,0)$ is a local maximum of the restriction of $\delta_\mathscr{C}$ to that  face   of $\mathscr{I}^n$.  The induction hypothesis  then   gives  $\delta_\mathscr{C}(1/2,\ldots,1/2,0)\le (n-3)/2$, and it follows that $\delta_\mathscr{C}(\vv{t})\le (n-2)/2$, a contradiction.   
\qed\end{potC2}

In view of the above,  from now on we  assume  that $p\ne 2$ and $n>2$.
Our proof will consist of four steps. 

 \subsubsection*{Proof outline:}
 
 \bigskip

 \begin{enumerate}
 \item   In Section \ref{ss: HKM and syzygy gaps} we  relate   $\delta^*_\mathscr{C}(\vv{c}/m)$, where $m<\sum_{i=1}^nc_i$, to the Hilbert--Kunz multiplicity of a   3-variable homogeneous polynomial $F=\ell^{\vv{c}}-z^mH$ with respect to the  ideal $\ideal{U,V,z}$, under the assumption that $\deg U=\deg V$ (or, equivalently, $\delta_\mathscr{C}(\vv{0})=0$).
  \item In Section \ref{ss: Brenner and Trivedi} we use results of Brenner and Trivedi to find another formula for that Hilbert--Kunz multiplicity,  thereby obtaining some information on  $\delta^*_\mathscr{C}(\vv{c}/m)$.

 \item In  Section \ref{ss: key lemma} we   prove that if $\delta_\mathscr{C}$ has a local maximum  at a point $\vv{a}/q$, where $\delta_\mathscr{C}(\vv{0})=0$, $\delta_\mathscr{C}(\vv{a}/q)\ge (n-1)/q$,  and  $d(\vv{a}/q,\vv{0})>1$, then   $\sum_{i=1}^na_i\equiv 0\pmod{p}$. That is done by choosing a convenient point $\vv{c}/m$, close enough to $\vv{a}/q$ to be ``under the effect'' of that local maximum, and  using the information on   $\delta^*_\mathscr{C}(\vv{c}/m)$ previously found.  Reflections then show that each corner $\vv{b}$   with  $\delta_\mathscr{C}(\vv{b})=0$ and $d(\vv{a}/q,\vv{b})>1$   yields  a congruence of the form  $\sum_{i=1}^n(\pm a_i)\equiv0\pmod{p}$. 
 \item We conclude  the proof in Section \ref{ss: conclusion}: assuming that  $\delta_\mathscr{C}$ has a local maximum at an interior  point $\vv{a}/p$, where it takes on a value  $\ge (n-1)/p$, we shall show that there are enough corners   $\vv{b}$ as above,    with  $\delta_\mathscr{C}(\vv{b})=0$ and $d(\vv{a}/p,\vv{b})>1$,  to guarantee that the corresponding congruences lead to  a contradiction.  Special cell classes are handled separately, through a simpler argument that takes advantage of their self-similarities.
 \end{enumerate}

\subsection{Hilbert--Kunz multiplicities and syzygy gaps}\label{ss: HKM and syzygy gaps}

Recall that  $\mathscr{C}=\cclass{U,V}$, where $U$ and $V$ have no common factor; in this subsection we add the extra assumption that $\deg U=\deg V=d$. We denote by   $\delta^*_\mathscr{C}$  the continuous extension of $\delta_\mathscr{C}$ to $[0,1]^n$, as in Remark \ref{extension by continuity}.

Let  $\vv{c}$ be a nonnegative integer vector, $G=\ell^\vv{c}$,  and $r=\deg G$. Fix $m$ with    $0<m<r$ such that $\vv{c}/m\in [0,1]^n$.   Let $H\in \Bbbk[x,y]$ be a homogeneous  polynomial of degree $r-m$,  prime to $G$,   and set  $F=G-z^mH\in \Bbbk[x,y,z]$.

\begin{defn} $\mu(F)$ is the Hilbert--Kunz multiplicity of   $\Bbbk[x,y,z]/\ideal{F}$ with respect to the ideal generated by the images of $U$, $V$, and $z$.
\end{defn}

We  shall  relate  $\mu(F)$ and  $\delta^*_\mathscr{C}(\vv{c}/m)$, proving:

\begin{thm}\label{HKM versus delta}
$\displaystyle \mu(F)=d r-\frac{r}{4}+ \frac{m^2}{4r}\cdot \delta_\mathscr{C}^*\left(\frac{\vv{c}}{m}\right)^2.$
\end{thm}

\begin{lem}\label{reduction to irreducible case}
Let $\lambda=\lambda(F)$ be the greatest divisor of $m$ for which $G/H$ is a $\lambda$th power in $\Bbbk(x,y)$. Then:
\begin{enumerate}
\item If $\lambda=1$, then $F$ is irreducible in $\Bbbk[x,y,z]$. 
\item If Theorem \ref{HKM versus delta} holds for $\lambda=1$, then it holds in general.
\end{enumerate}
\end{lem}

\begin{pf}
 Since $G$ and $H$ are relatively prime, any nontrivial factorization of  $F=G-z^mH$   in $\Bbbk[x,y,z]$ would have factors of degree $<m$ in $z$, giving a nontrivial factorization of $z^m-G/H$  in $\Bbbk(x,y)[z]$. But   $z^m-G/H$ is irreducible in $\Bbbk(x,y)[z]$ if $\lambda=1$ (see, \eg \cite[Chapter VI, Theorem 9.1]{lang});  this gives us 1.

Suppose now $\lambda>1$. Since $G$ and $H$ are relatively prime, both $G$ and $H$ are $\lambda$th powers, and we can write $z^\lambda-G/H=\prod_{i=1}^\lambda(z-\ell^{\vv{c}/\lambda}/H_i)$, where the $H_i$  are  $\lambda$th roots of $H$. Replacing $z$ with $z^{m/\lambda}$ and multiplying through by $H$ we  see that  $F=F_1\cdots F_\lambda$, where $F_i=\ell^{\vv{c}/\lambda}-z^{m/\lambda}H_i$. If Theorem \ref{HKM versus delta} holds for $\lambda=1$, it   gives formulas for the Hilbert--Kunz multiplicity $\mu(F_i)$ of each $F_i$. But the additivity of the Hilbert--Kunz multiplicity shows that $\mu(F)=\sum_{i=1}^\lambda\mu(F_i) $, and  adding up those formulas gives Theorem \ref{HKM versus delta} for $F$. 
\qed
\end{pf}

In the remainder of this section we assume that $\lambda(F)=1$, so  $F$ is irreducible in $\Bbbk[x,y,z]$.
Let $\bar{z}$ and $w$ be elements in an extension of $\Bbbk(x,y)$ such that $\bar{z}^m=G$ and $w^m=H$.  Set $\bar{x}=wx$ and $\bar{y}=wy$. Then 
\begin{align*}
F(\bar{x},\bar{y},\bar{z})&=G(\bar{x},\bar{y})-\bar{z}^mH(\bar{x},\bar{y})\\
&=w^r G(x,y)-G(x,y)\cdot w^{r-m}H(x,y)\\
&=w^r G(x,y)-G(x,y)\cdot w^{r-m}\cdot w^m\\
&=0, 
\end{align*}
so $\Bbbk[\bar{x},\bar{y},\bar{z}]$ is a homogeneous coordinate ring for $F$. 
We now consider the rings in the following diagram.
\begin{center}
\bigskip
\unitlength .8 mm
\begin{picture}(40,35)(0,0)
\linethickness{0.3mm}
\multiput(0,30)(0.15,-0.24){83}{\line(0,-1){0.24}}
\multiput(20,10)(0.15,0.24){83}{\line(0,1){0.24}}
\put(0,35){\makebox(0,0)[cc]{$A=\Bbbk[x,y]$}}

\put(30,35){\makebox(0,0)[cc]{$B=\Bbbk[\bar{x},\bar{y},\bar{z}]$}}

\put(46.5,2){\makebox(0,0)[cc]{$\begin{array}{rcl}R&=&\Bbbk[U(\bar{x},\bar{y})^m,V(\bar{x},\bar{y})^m,\bar{z}^{m}]\\&=&\Bbbk[H^dU^m,H^dV^m,G]\end{array}$}}

\put(29,20){\makebox(0,0)[cc]{$\beta$}}

\put(3,20){\makebox(0,0)[cc]{$\alpha$}}

\end{picture}
\bigskip
\end{center}
Here $\alpha$ and $\beta$ denote the ranks of $A$ and $B$ over $R$; these are finite, according to the next lemma. 

\begin{lem}
$A$ and $B$ are finite over $R$. 
\end{lem}

\begin{pf}
The ideal $\ideal{H^dU^m,H^dV^m,G}$  of $A$ is $\ideal{x,y}$-primary, as its generators have no common factor; so it contains $x^s$ and $y^s$ for some $s$. Let $M$ be the $R$-submodule of $A$ generated by   $x^iy^j$, with $i+j\le 2s-2$.  Writing    $x^s$ and  $y^s$ as $A$-linear combinations of  $H^dU^m$, $H^dV^m$, and $G$, we see that  they can be expressed  as $R$-linear combinations of monomials of degree $<s$.  It follows easily  that $xM\subseteq M$ and $yM\subseteq M$,  so that any monomial in $x$ and $y$ is in $M$. Hence $A=M$, and $A$ is finite over $R$.

Arguing along the same lines,  choosing $s$ such that $x^s$ and $y^s$ are in the ideal
   $\ideal{U^m,V^m}$ of $A$ we can show that $B$ is generated over $R$ by  monomials  $\bar{x}^i\bar{y}^j\bar{z}^k$, with $i+j\le 2s-2$ and $k<m$.  
   \qed
\end{pf}

\begin{defn}
For any nonnegative integer $k$, $\mu(k)$ is the Hilbert--Kunz multiplicity of $B$ with respect to the ideal  $J(k)=\ideal{U(\bar{x},\bar{y})^k,V(\bar{x},\bar{y})^k,\bar{z}^{k}}$. 
\end{defn}

\begin{lem}\label{HKM versus deg J} $\mu(km)=(\beta/\alpha)\cdot\deg J$, where  $J$ is  the ideal of $A$ generated by $H^{kd}U^{km}$, $H^{kd}V^{km}$, and $G^{k}$. 
\end{lem}

\begin{pf}
The generators of $J$ are precisely the generators of $J(km)$, and are elements of $R$; let $I$ be the ideal they generate in $R$. Then $\mu(km)$ coincides with the Hilbert--Kunz multiplicity  of $B$ (seen as an $R$-module)  with respect to $I$. Using Theorem 1.8 of  \cite{paul1} we see that this Hilbert--Kunz  multiplicity is just   $\beta/\alpha$ times  the Hilbert--Kunz multiplicity of $A$ with
respect to $I$. But this is $(\beta/\alpha) \cdot \deg IA$, since $A=\Bbbk[x,y].$
\qed
\end{pf}

\begin{lem}\label{beta/alpha}
$\beta/\alpha=m^2/r$.
\end{lem}

\begin{pf} 
Let \[K=\Bbbk\left(\frac{U(\bar{x},\bar{y})^m}{V(\bar{x},\bar{y})^m},\frac{V(\bar{x},\bar{y})^m}{\bar{z}^{dm}}\right)=\Bbbk\left(\frac{U^m}{V^m},\frac{H^dV^m}{G^{d}}\right)\subseteq \Bbbk\left(\frac{x}{y}\right)= \Bbbk\left(\frac{\bar{x}}{\bar{y}}\right).\]

The field of fractions of $R$ is $K(G)=K(\bar{z}^m)$, and $\alpha$ and $\beta$ are the degrees    of $\Bbbk(x,y)$ and $\Bbbk(\bar{x},\bar{y},\bar{z})$ over that field.
Because    $\bar{z}$ is a root of the degree $m$ irreducible polynomial $F(\bar{x},\bar{y},z)\in \Bbbk(\bar{x},\bar{y})[z]$ (see the proof of Lemma \ref{reduction to irreducible case}),  we have
\begin{align*}
\beta&=\left[\Bbbk\left(\bar{x},\bar{y},\bar{z}\right): \Bbbk\left(\bar{x},\bar{y}\right)\right]\cdot  \left[\Bbbk\left(\bar{x},\bar{y}\right):\Bbbk\left(\bar{x}/\bar{y},\bar{y}^m\right)\right]\cdot \left[\Bbbk\left(\bar{x}/\bar{y},\bar{y}^m\right):K\left(G\right)\right]\\
&=m^2\cdot \left[\Bbbk\left(\bar{x}/\bar{y},\bar{y}^m\right):K\left(G\right)\right].
\end{align*}
But $\bar{y}^m/G=w^my^m/G=Hy^m/G\in \Bbbk(x/y)=\Bbbk(\bar{x}/\bar{y})$, so $\Bbbk(\bar{x}/\bar{y},\bar{y}^m)=\Bbbk(x/y,G)$, and 
\begin{equation}\label{beta}
\beta=m^2\cdot \left[\Bbbk\left(x/y,G\right):K(G)\right].
\end{equation}
A similar calculation gives
\begin{align}\label{alpha}
\alpha &= \left[\Bbbk(x,y):\Bbbk\left(x/y,y^r\right)\right]\cdot\left[\Bbbk\left(x/y,y^r\right):K(G)\right]\nonumber\\
&=r\cdot \left[\Bbbk\left(x/y,G\right):K(G)\right],
\end{align}
and comparing (\ref{beta}) and (\ref{alpha}) we get the desired result.
\qed
\end{pf}

\begin{cor}\label{discrete formula}
$\displaystyle \mu(km)=k^2m^2dr-\frac{k^2m^2r}{4}+\frac{m^2}{4r}\cdot \delta(U^{km},V^{km},G^k)^2.$
\end{cor}

\begin{pf} Note that $Q(kdr,kdr,kr)=4 d k^2 r^2-k^2 r^2$, where $Q$ is the quadratic form of Proposition \ref{formula1}. So  
Lemmas \ref{HKM versus deg J} and \ref{beta/alpha}, together with Proposition \ref{formula1}, give us 
\[\mu(km)=\frac{m^2}{4r}\left( 4 d k^2 r^2-k^2 r^2 +\delta(H^{kd}U^{km},H^{kd}V^{km},G^{k})^2\right).\]
But    $\delta(H^{kd}U^{km},H^{kd}V^{km},G^{k})=\delta(U^{km},V^{km},G^{k})$, by Proposition \ref{deltaproperties1}.
\qed
\end{pf}

We now need a ``continuous version'' of the above result; we shall arrive at the desired formula for $\mu(F)$ by replacing $k$ with $1/m$ in that continuous version. 

\begin{defn}
$\delta_G^*:[0,1]^3\to \mathbb{R}$ is the continuous  function  such that
\[\delta_G^*\left(\frac{\vv{a}}{q}\right)=\frac{1}{q}\cdot \delta(U^{a_1},V^{a_2},G^{a_3}),\]
for any $q$ and any   $\vv{a}\in [q]^3$.
\end{defn}

Directly from the definition of Hilbert--Kunz multiplicity it follows that $\mu(pk)=p^2\cdot \mu(k)$, so we may define  a function $\mathscr{I}\to \mathbb{Q}$, $k/q\mapsto q^{-2}\cdot \mu(k)$. This function is  uniformly  continuous (see  \ref{appendix} for a proof in a more general setting), so we can extend it to a continuous  function on $[0,1]$.

\begin{defn}
$\mu^*:[0,1]\to \mathbb{R}$ is the continuous   extension of the function $k/q\mapsto q^{-2}\cdot \mu(k)$. 
\end{defn}

\begin{cor}\label{continuous formula}
$\displaystyle \mu^*(tm)=t^2m^2dr-\frac{t^2m^2r}{4}+\frac{m^2}{4r}\cdot \delta_G^*(tm,tm,t)^2$, for all $t\in [0,1/m].$
\end{cor}

\begin{pf}
Corollary \ref{discrete formula} gives the formula for  $t\in [0,1/m]\cap\mathbb{Z}[1/p]$, and the result follows by  continuity.
\qed
\end{pf}

We can now complete the proof of Theorem \ref{HKM versus delta}. 

\begin{proofHKMformula}
  Since,  for any $q$ and   $a\in [q]$,
  \begin{align*}
  \delta_G^*\left(1,1,\frac{a}{q}\right)&=\frac{1}{q}\cdot \delta(U^q,V^q,G^a)\\
  &=\frac{1}{q}\cdot \delta(U^q,V^q,\ell^{a\vv{c}})\\
  &=\delta_\mathscr{C}^*\left(\frac{a}{q}\cdot \vv{c}\right),
  \end{align*}
       the continuity of  $\delta_G^*$ and $\delta_\mathscr{C}^*$ implies that  $\delta_G^*(1,1,t)=\delta_\mathscr{C}^*(t\vv{c})$, for any   $t\in [0,1]$.  Setting $t=1/m$ in the identity of Corollary \ref{continuous formula} we get the desired result:
\[
\mu(F)=\mu(1)=dr-\frac{r}{4}+\frac{m^2}{4r}\cdot \delta_G^*\left(1,1,\frac{1}{m}\right)^2
=dr-\frac{r}{4}+\frac{m^2}{4r}\cdot \delta_\mathscr{C}^*\left(\frac{\vv{c}}{m}\right)^2.\]
\qed
\end{proofHKMformula}

\subsection{An application of sheaf theory}\label{ss: Brenner and Trivedi}
 
 \begin{defn} 
Let $F\in \Bbbk[x,y,z]$ be an irreducible homogeneous polynomial, and $\projcurve$ be a desingularization of the projective curve defined by $F$. Then $\gamma(F)=2\genus(\projcurve)-2$.
\end{defn}
 
The following result will be essential to our argument:

\begin{thm}\label{brenner-trivedi}  Let  $q>1$ be a power of $p$.
Let $F\in \Bbbk[x,y,z]$ be an irreducible degree $r$ homogeneous polynomial, and let $\mu(F)$ be the Hilbert--Kunz multiplicity of $\Bbbk[x,y,z]/\ideal{F}$ with respect to a zero-dimensional ideal  $I$ generated by three homogeneous elements of degrees $d_1$, $d_2$, and $d_3$. Then 
\begin{equation}\label{conclusion}
\mu(F)=\frac{r}{4}\cdot Q(d_1,d_2,d_3)+\frac{l^2}{4r},
\end{equation}
where  $l$ is a number in  $\mathbb{Z}[1/p]$ such that  $ql\in p\mathbb{Z}$ or $0<ql\le \gamma(F)$, and $Q$ is the quadratic form   of Proposition \ref{formula1}, 
\[Q(d_1,d_2,d_3)  =
2d_1d_2+2d_1d_3+2d_2d_3-{d_1}^2-{d_2}^2-{d_3}^2.\] 
\end{thm}

When $d_1=d_2=d_3=1$ this is Theorem 5.3 of Trivedi \cite{trivedi}. The general case is treated similarly, but now we need a  result from Brenner \cite{brenner} and a lemma of Trivedi. 
  Before we give the proof, we recall some of the terminology used in those papers. For a rank $r$ vector bundle $\vb{S}$ on a smooth projective curve  $\projcurve$ over an algebraically closed field, $\deg(\vb{S})$ is the degree of the line bundle $\bigwedge^r\vb{S}$; the degree is additive in the category of vector bundles on $\projcurve$. The \emph{slope} of $\vb{S}$  is defined as $\deg(\vb{S})/r$. The vector bundle $\vb{S}$ is \emph{semistable} if $\slope(\vb{T})\le \slope(\vb{S})$, for every subbundle $\vb{T}$ 
of $\vb{S}$.  $\vb{S}$ is \emph{strongly semistable} if  its pull-back by each $e$th iterate  of the absolute Frobenius $\frob:\projcurve\to \projcurve$ is semistable. 

\begin{pf}
 Let $B=\Bbbk[x,y,z]/\ideal{F}$, and let  $R$ be the integral closure of $B$. The  Hilbert--Kunz multiplicities of $B$
    and   $R$ with respect to $I$ are equal, and $\projcurve= \proj R$ is the desingularization of the  projective  curve defined
    by $F$.  
 
 In \cite[Corollary 4.4]{brenner} Brenner considers a rank 2 vector bundle $\vb{S}$ on $\projcurve$---the pull-back to $\projcurve$  of the bundle of syzygies between the three
                               homogeneous generators of $I$. The degree of $\vb{S}$  is  $-(d_1+d_2+d_3)r$. He shows that if $\vb{S}$ is strongly semistable, then \eqref{conclusion} holds with $l=0$.\footnote{Brenner makes the assumption that $R$ is generated by finitely many elements of degree 1, which is not necessarily the case  here,  but that assumption can be weakened---that is the content of his footnote 1.}  If, on the other hand, $\vb{S}$ is not strongly semistable,  let  $e$  be the least number for which $\frob^{e*}(\vb{S})$ is not semistable. Then $\frob^{e*}(\vb{S})$ has a subbundle $\vb{L}$  with $\slope(\vb{L}) >\slope(\frob^{e*}(\vb{S}))$.  Because   $\vb{S}$  has rank  2,  $\vb{L}$ and $\vb{M}=\frob^{e*}(\vb{S})/\vb{L}$ are line bundles, and the condition on the slopes is equivalent to $\deg(\vb{L})>(\deg(\vb{L})+\deg(\vb{M}))/2$, or  $\deg(\vb{L})>\deg(\vb{M})$.

 Brenner  then   sets 
\[\nu_1=-\frac{\deg(\vb{L})}{rq^*}\ \ \ \ \text{and}\ \ \ \ \nu_2=-\frac{\deg(\vb{M})}{rq^*},\]
 where $q^*=p^e$, and shows that 
\begin{equation}\label{brenner result}
\mu(F)=r\left(\nu_2^2-\nu_2\sum_{i=1}^3d_i+\sum_{i<j}d_id_j\right).
\end{equation}
Note that    $\deg(\vb{L})+\deg(\vb{M})=\deg (\frob^{e*}(\vb{S}))=-(d_1+d_2+d_3)rq^*$,  so 
   $\nu_1+\nu_2=d_1+d_2+d_3$.  Using this, Eq. \eqref{brenner result} gives us \eqref{conclusion} with $l=r(\nu_2-\nu_1)=(\deg(\vb{L})-\deg(\vb{M}))/q^*$.  
 
  If $q>q^*$, then $ql\in p\mathbb{Z}$, and we are done. If  $q\le q^*$,   Lemma 5.2 of Trivedi \cite{trivedi} comes into play. Since $e$ was chosen  to be  the least number for which  $\frob^{e*}(\vb{S})$ is not semistable, Trivedi's result says that $\deg(\vb{L})-\deg(\vb{M})\le \gamma(F)$. Since $q\le q^*$, it follows that $ql=q(\deg(\vb{L})-\deg(\vb{M}))/q^*\le \gamma(F)$.
\qed\end{pf}

 Now let $\mathscr{C}$, $U$, $V$, $d$, $\vv{c}$, $r$, and $m$  be as at the start of Section \ref{ss: HKM and syzygy gaps}.  Comparing the above theorem to Theorem \ref{HKM versus delta}  we shall obtain some information on $\delta_\mathscr{C}^*(\vv{c}/m)$ which will play an important role in the next section. 
 
  For ease of notation, for any vector $\vv{a}=(a_1,\ldots,a_n)$ we  write $\norm{a}=\sum_{i=1}^n|a_i|$; we shall refer to  $\norm{a}$ as the \emph{norm} of the vector $\vv{a}$.

\begin{lem}\label{hurwitz}
  Suppose that some $k>1$ divides each $c_i$; write $\vv{c}=k\vv{a}$. Suppose further that $m$ is prime to $p$ and to $k$, and divisible  by  $\norm{a}$. Let $q>1$ be a power of $p$. Then one of the following holds:

\begin{enumerate}
\item  $\displaystyle{qm\cdot \delta_\mathscr{C}^*\left(\frac{\vv{c}}{m}\right)\in p\mathbb{Z}}$
\item $\displaystyle{\delta_\mathscr{C}^* \left(\frac{\vv{c}}{m}\right) <\frac{n-1}{q} - \frac{\norm{a}}{qm}}$
\end{enumerate}
\end{lem}

\begin{pf}   Let  $\lambda$ be the greatest common divisor  of $m$ and the $c_i$. Since $m$ is prime to $k$, so is $\lambda$, and $\lambda$ divides each $a_i$. If we replace  $\vv{c}$, $m$, and $\vv{a}$ by their
        quotients by $\lambda$,  then $\vv{c}/m$ and $\norm{a}/m$ are unchanged. So it suffices to show that 1 or 2 holds after this replacement, and we may assume         that  $\lambda=1$. Now set  $F=\ell^\vv{c}-z^mL^{r-m}$, where $L\in \Bbbk[x,y]$ is a linear form prime to each $\ell_i$. $F$ is an  irreducible homogeneous polynomial of degree $r$ (irreducibility follows from Lemma \ref{reduction to irreducible case}, since $\lambda=1$).

 We now apply Theorem  \ref{brenner-trivedi} with $I=\ideal{U, V, z}$, to find that $\mu(F) =d r-r/4+l^2/(4r)$,    where $ql\in p\mathbb{Z}$ or $ 0<ql \le  \gamma(F)$.  Comparing with     Theorem \ref{HKM versus delta}  we see that $l=m\cdot \delta_\mathscr{C}^*\left(\vv{c}/m\right)$.  So either $qm\cdot \delta_\mathscr{C}^*(\vv{c}/m)\in p\mathbb{Z}$ or $\delta_\mathscr{C}^*(\vv{c}/{m}) \le  \gamma(F)/(qm)$.
It only remains to show that $\gamma(F) <(n-1)m  - \norm{a}$.

  The desingularization $\projcurve$ of the projective curve defined by $F$ is an $m$-sheeted branched covering of $\mathbb{P}^1$, tamely ramified, since $m$ is prime to $p$. According to  the Hurwitz formula, \[\gamma(F)=-2m+(\text{terms coming from ramification}).\]  Ramification can only occur at zeros of the $\ell_i$ and of $L$.  Because of tameness, the contribution from the zero of each  $\ell_i$ is at most $m-1$. Now note that the greatest common divisor of $m$ and $r-m$ is $\norm{a}$, since $r/\norm{a}=k$, while $m/\norm{a}$ is an integer prime to $k$. So
   over the zero of $L$ there are $\norm{a}$ points of $\projcurve$, each of ramification degree $m/\norm{a}$, providing a  contribution of $m-\norm{a}$ to $\gamma(F)$. So
\[\gamma(F)\le -2m+n(m-1)+m-\norm{a}<(n-1)m-\norm{a}.\]  
\qed\end{pf}

\subsection{A key lemma}\label{ss: key lemma}

The following lemma  will play a crucial  role in our proof of Theorem \ref{bounds1}. 

\begin{lem}\label{main lemma}
Suppose $\delta_\mathscr{C}(\vv{0})=0$ and $\delta_\mathscr{C}$ has a local maximum at $\vv{a}/q$, where  $\delta_\mathscr{C}(\vv{a}/q)\ge (n-1)/q$  and $\norm{a}> q>1$. Then $p$ divides $\norm{a}$. 
\end{lem}

\begin{pf}
Suppose not.  As discussed in Section \ref{reductions}, an inductive argument allows us to assume that $\vv{a}/q$ is an interior point of $\mathscr{I}^n$, \ie $0<a_i<q$, for all $i$. Note that $\deg U=\deg V$, since $\delta_\mathscr{C}(\vv{0})=0$, so  we are in the situation of  Section  \ref{ss: HKM and syzygy gaps}.  By looking at values of $\delta^*_\mathscr{C}$ at   conveniently chosen points $\vv{c}/m$ that are sufficiently close to $\vv{a}/q$ to be ``under the influence'' of that local maximum (see Theorem \ref{max-corner}) and using  Lemma  \ref{hurwitz}, we  shall  prove that $\Delta(\vv{a}):=q\cdot\delta_\mathscr{C}(\vv{a}/q)$ is congruent modulo ${p}$ to both $\norm{a}$ and $-\norm{a}$. This will give us a contradiction, since we are assuming that $p\ne 2$.

Since $p$ does not divide $\norm{a}$, we can find a multiple $m$ of  $\norm{a}$  of the form $m=kq+1$, with $k>1$.  Note that our assumptions on $\vv{a}$ imply that $ka_i<m<k \norm{a}$. Let $\vv{c}=k \vv{a}$; then 
\[\frac{\vv{a}}{q}-\frac{\vv{c}}{m}=\frac{m\vv{a}-kq\vv{a}}{qm}=\frac{\vv{a}}{qm},\]
so
\[d\left(\frac{\vv{a}}{q},\frac{\vv{c}}{m}\right)=\frac{\norm{a}}{qm}. \]
But since $k a_i<m$, it follows that $k \norm{a}<nm$, so $\norm{a}/m<n/k\le n/2\le n-1$. Thus
\[d\left(\frac{\vv{a}}{q},\frac{\vv{c}}{m}\right)<\frac{n-1}{q}\le \delta^*_\mathscr{C}\left(\frac{\vv{a}}{q}\right),\]
and Theorem \ref{max-corner} (and continuity) shows that 
\begin{equation}\label{eq5}
\delta^*_\mathscr{C}\left(\frac{\vv{c}}{m}\right)=\delta^*_\mathscr{C}\left(\frac{\vv{a}}{q}\right)-d\left(\frac{\vv{a}}{q},\frac{\vv{c}}{m}\right)
=\delta^*_\mathscr{C}\left(\frac{\vv{a}}{q}\right)-\frac{\norm{a}}{qm}.
\end{equation}

 Since $\delta_\mathscr{C}^* (\vv{a}/q) \ge (n-1)/q$,  situation 2 of  Lemma \ref{hurwitz} cannot hold, and so $qm\cdot \delta_\mathscr{C}^* (\vv{c}/m)\in p\mathbb{Z}$. Multiplying \eqref{eq5} through by
       $qm$ we find that  $m \Delta(\vv{a})-\norm{a} \in  p\mathbb{Z}$. Since $m\equiv 1 \pmod{p}$, $\Delta(\vv{a})\equiv \norm{a} \pmod{p}$. By repeating the argument with an $m$ that is 
       divisible by  $\norm{a}$ and of the form $kq-1$ we find that $\Delta(\vv{a})\equiv -\norm{a} \pmod{p}$. So $p$ divides $\norm{a}$, contradicting our assumption.      
\qed\end{pf}

Using reflections, the following corollary is immediate from Lemma \ref{main lemma}.

\begin{cor}\label{cor to main lemma}
Let $\vv{c}= (\epsilon_1,\ldots,\epsilon_n)$ be a corner of $\mathscr{I}^n$. Suppose $\delta_\mathscr{C}(\vv{c})=0$ and $\delta_\mathscr{C}$ has a local maximum at $\vv{a}/q$, where $q>1$, $\delta_\mathscr{C}(\vv{a}/q)\ge (n-1)/q$,  and  $d(\vv{a}/q,\vv{c})> 1$. Then $p$ divides $\sum_{i=1}^n(-1)^{\epsilon_i}a_i$. 
\end{cor}

\begin{rem}\label{rem: assumption on distance may be dropped} In the case of a special cell class $\mathscr{C}$, self-similarity properties allow us to drop the assumption that $\norm{\vv{a}}>q$ in Lemma \ref{main lemma}. In fact, suppose $\delta_\mathscr{C}(\vv{0})=0$ (so $\mathscr{C}=\cclass{x,y}$, by Proposition \ref{special cells are unique}) and   $\delta_\mathscr{C}$ has a local maximum at $\vv{a}/q$, where where $q>1$ and   $\delta_\mathscr{C}(\vv{a}/q)\ge (n-1)/q$.  Setting $\vv{b}=(p-1,p-1,0,\ldots,0)$,
  calculations made in Example  \ref{self-sim 1} show that $\mathscr{C}=T_{p|\vv{b}}\mathscr{C}$. So \[\delta_\mathscr{C}(\vv{t})=(T_{p|\vv{b}}\delta_\mathscr{C})(\vv{t})=p\cdot \delta_\mathscr{C}\left(\frac{\vv{t}+\vv{b}}{p}\right),\]  for all $\vv{t}\in \mathscr{I}^n$.  So $\delta_\mathscr{C}$ also has a local maximum at $(\vv{a}+q\vv{b})/(pq)$, where it takes on a value $\ge (n-1)/(pq)$. Setting $\vv{a^*}=\vv{a}+q\vv{b}$, we see that $\vv{a^*}/(pq)$ satisfies the hypotheses of Lemma \ref{main lemma}, so  $p$ divides $\norm{\vv{a^*}}$; but $\norm{\vv{a}}\equiv \norm{\vv{a^*}}\pmod{p}$.
\end{rem}

\subsection{Concluding the proof of Theorem \ref{bounds1}}\label{ss: conclusion}

We have now the machinery necessary  to prove Theorem \ref{bounds1}, which we restate below: 

\begin{thmCpreview}
Suppose $\delta_\mathscr{C}$ has a local maximum at $\vv{a}/q$, where $q>1$ and $\vv{a}/q$ is reduced, in the sense  that some coordinate  $a_i/q$ is  reduced.  Then    \[\delta_\mathscr{C}\left(\frac{\vv{a}}{q}\right)\le \frac{n-2}{q}.\] 
\end{thmCpreview}

We start by considering the particular case of special cell classes. As observed in Remark \ref{rem: bound for special cell classes}, in this case Theorem \ref{bounds1} is   equivalent to  Monsky's result from \cite{mason}. However,  with the machinery   already developed  its proof is simple enough, so we include it here. (Remark \ref{rem: assumption on distance may be dropped} and further self-similarity properties make this special case a lot less convoluted than the general case.)

\begin{potC3} 
In view of  Proposition \ref{special cells are unique} we may  assume   $\mathscr{C}=\cclass{x,y}$. Suppose $\delta_\mathscr{C}$ has a local maximum at $\vv{a}/q$, with $q>1$ and $\delta_\mathscr{C}(\vv{a}/q)\ge (n-1)/q$. Remark \ref{rem: assumption on distance may be dropped} allows us to use Lemma \ref{main lemma} to conclude that 
\begin{equation}\label{first cong}
a_1+a_2+\cdots+a_n=\norm{\vv{a}}\equiv 0\pmod{p}.  
\end{equation}

Now choose  $q^*\ge n-2$ and set $\vv{b}=(q^*-1,q^*-n+2,0,\ldots,0)$; Eq. (\ref{eqn: self-sim}) of Example \ref{self-sim 1} shows that 
$\delta_\mathscr{C}=(T_{q^*|\vv{b}}\delta_\mathscr{C})\circ R_2\cdots R_n,$
so 
\[\delta_\mathscr{C}(\vv{t})=q^*\cdot \delta_\mathscr{C}\left(\frac{R_2\cdots R_n(\vv{t})+\vv{b}}{q^*}\right),\]
for all  $\vv{t}\in \mathscr{I}^n$.
In particular,   $\delta_\mathscr{C}$ also has a local maximum at  $(R_2\cdots R_n(\vv{a}/q)+\vv{b})/q^*=(qR_2\cdots R_n(\vv{a}/q)+q\vv{b})/(qq^*)$, where it takes on a value $\ge (n-1)/(qq^*)$.
Lemma \ref{main lemma} then shows that 
\begin{equation}\label{second cong}
a_1-a_2-\cdots-a_n\equiv \|qR_2\cdots R_n(\vv{a}/q)+q\vv{b}\|\equiv 0\pmod{p}.
\end{equation}

Combining (\ref{first cong})  and  (\ref{second cong}) we conclude that $2a_1\equiv 0\pmod{p}$, and since we are assuming that $p\ne 2$, $a_1\equiv 0\pmod{p}$. Similarly, we show that $p$ divides each of the other $a_i$, so $\vv{a}/q$ is not reduced.
\qed\end{potC3}

We now turn to the proof of  Theorem \ref{bounds1} for arbitrary cell classes. The following simple lemmas will be helpful in our argument. 

\begin{lem}\label{PL or vanishes at another layer} Suppose $\delta_\mathscr{C}$ vanishes at all corners of norm $k$, for some $k$. Then one of the following holds:
\begin{enumerate}
\item $\delta_\mathscr{C}$ also vanishes at some corner of norm $k+2$ or $k-2$.
\item $\delta_\mathscr{C}$ is piecewise linear, with local maxima only at the origin $\vv{0}$ and at its opposite corner, $\vv{1}$. 
\end{enumerate} 
\end{lem}

\begin{pf}  Since  $\delta_\mathscr{C}=0$  at all   corners of norm $k$, $\delta_\mathscr{C}=1$ at all corners of norm $k\pm1$.
If we are not in situation 1, then  $\delta_\mathscr{C}=2$ at all corners of norm $k\pm 2$. Then Proposition \ref{2Dconv} forces $\delta_\mathscr{C}$  to be $3$ at all corners of norm $k\pm 3$,  $4$ at all corners of norm $k\pm 4$,  and so on.
Thus $\delta_\mathscr{C}|_{\Xset_1}$ has local maxima at $\vv{0}$ and $\vv{1}$, where it takes on the values $k$ and $n-k$. By Theorem \ref{max-corner}, the same is true for $\delta_\mathscr{C}$, and   $\delta_\mathscr{C}(\vv{t})=\max 
\{k-d(\vv{t}, \vv{0}), n-k-d(\vv{t},\vv{1})\}.$ 
\qed
\end{pf}

\begin{lem}\label{congruences} Suppose $\delta_\mathscr{C}$ has a local maximum at  an interior point $\vv{a}/p$ of $\Xset_p$, where $\delta_\mathscr{C}(\vv{a}/p)\ge (n-1)/p$.
Furthermore, suppose $\delta_\mathscr{C}$ vanishes at all corners of norm $k$, for some $k$ with $2\le k\le n-2$.   Then all the $a_i$ are congruent modulo $p$, and  
\[(n-2k)a_i\equiv 0\pmod{p}.\]
The same conclusion   holds  if $k=1$, provided the distance from each of the corners of norm $1$ to $\vv{a}/p$ is $>1$.
\end{lem}

\begin{pf}
The local maximum $\vv{a}/p$ can be  within distance $1$ of at most one of the corners of norm $k$;   Corollary \ref{cor to main lemma} gives a   linear congruence modulo  $p$ for each of the $\binom{n}{k}$ or $\binom{n}{k}-1$ corners of norm $k$ that are ``far'' from $\vv{a}/p$. For each $i\ne j$ there are two   congruences that differ only by the signs of $a_i$ and $a_j$ (more precisely, there are   $\binom{n-2}{k-1}$ or $\binom{n-2}{k-1}-1$ such pairs).  Subtracting one such congruence  from the other we find that $2(a_i-a_j)\equiv 0\pmod{p}$, and since $p\ne 2$,   $a_i\equiv  a_j\pmod{p}$.     
Substituting that into any of the congruences we find $k(-a_i)+(n-k)a_i\equiv 0\pmod{p}$,  giving the result.  If $k=1$, the same argument applies, but we need congruences associated to \emph{all} corners of norm 1, hence the need for the extra assumption.
\qed
\end{pf}

We   can now conclude the proof of Theorem \ref{bounds1}.     

\begin{potC}    As pointed out in Section \ref{reductions}, we may assume  $n\ge 3$ and   $q=p$, and an inductive argument allows us to  restrict our attention to  interior points.  
Aiming at a contradiction, suppose $\delta_\mathscr{C}$ has a local maximum  at  an interior point $\vv{a}/p$, where  it takes on a value $\ge (n-1)/p$.

We would like to  arrange  a situation where we can use Lemma \ref{congruences}. By using a reflection, which changes the $a_i$  (modulo $p$) only by a sign, we may assume that the restriction of $\delta_\mathscr{C}$ to the corners of $\mathscr{I}^n$ attains its maximum value at the origin; let $k=\delta_\mathscr{C}(\vv{0})$.  
Note that if  $k\ge n$, then
 $\delta_\mathscr{C}$ is linear and its only local maximum is at  the origin,  while if $k=n-1$ then  $\delta_\mathscr{C}$ is piecewise  linear with local maxima only at  the origin and its opposite corner. In either case, 
 the existence of the local maximum at $\vv{a}/p$ is contradicted; so  henceforth  we assume $k\le n-2$. Theorem \ref{max-corner} then shows that $\delta_\mathscr{C}$ vanishes at all corners of norm $k$.

Difficulties may arise if $k=1$,  as  Lemma \ref{congruences} would then require the distance between $\vv{a}/p$ and each corner of norm 1 to be $>1$.  But these difficulties may be dealt with by using further reflections. Suppose, for instance, that $d(\vv{a}/p,(1,0,\ldots,0))\le 1$. Since the maximum value that $\delta_\mathscr{C}$ takes on at the corners is 1,   $\delta_\mathscr{C}$   vanishes at all   corners with an odd norm. The distance between $\vv{a}/p$ and each  such corner other than $(1,0,\ldots,0)$ is $>1$. Replacing $\mathscr{C}$ with  $R_2 R_3\mathscr{C}$ we arrive at   the desired situation: $\delta_\mathscr{C}$ now vanishes at all corners of norm 1, and    the distance between $\vv{a}/p$ and each of these corners is  $>1$. Lemma \ref{congruences} can thus be used even  if  $k=1$. (Note that    what made it possible for us to get around the difficulties    was the existence  of     an extra  ``layer'' of zeros of $\delta_\mathscr{C}$, namely the corners of norm 3.)
 
 Applying Lemma \ref{congruences}  we obtain
\begin{equation}\label{cong2}
(n-2k)a_i\equiv 0\pmod{p}.
\end{equation}
Since situation 2 of Lemma \ref{PL or vanishes at another layer} contradicts the existence of the local maximum at $\vv{a}/p$,  we may assume that $\delta_\mathscr{C}$ also vanishes at a corner of norm $k+2$. If the distance from $\vv{a}/p$ to that corner is $>1$, Corollary \ref{cor to main lemma} gives us another congruence $(n-2k-4)a_i\equiv 0\pmod{p}$; together with (\ref{cong2}), this shows that  $p$ divides $a_i$, a contradiction. 
If the distance  between $\vv{a}/p$ and that corner is $\le 1$, then $p\ge n-1$, since that distance is at least $ \delta_\mathscr{C}(\vv{a}/p)$, and  $\delta_\mathscr{C}(\vv{a}/p)\ge (n-1)/p$. So $p$ cannot divide $n-2k$, and (\ref{cong2}) gives us a contradiction, \emph{unless} $k=n/2$, in which case (\ref{cong2}) is of no help.

It remains to deal with the case $k=n/2$. In this case,    among all corners with norm $\ge n/2$ we choose a corner $\vv{c}$  where $\delta_\mathscr{C}$ is maximum; let $k'=\delta_\mathscr{C}(\vv{c})$. We may assume  that $k'<n/2$, as otherwise we would be in situation 2 of Lemma \ref{PL or vanishes at another layer}. If $k'>1$,  we use a  reflection to bring  $\vv{c}$ to  the origin, and conclude the proof by arguing exactly as above. If $k'=1$ we do the same,  with  some extra care---we need to choose $\vv{c}$ with $\norm{\vv{c}}\ge n/2+3$.  This ensures that all corners at distance 1 or 3 from $\vv{c}$ have norm $\ge n/2$, so   that $\delta_\mathscr{C}$ vanishes at  all those corners, guaranteeing  that extra ``layer'' of zeros   needed  for the workaround in the third paragraph of the proof.    Finding a  $\vv{c}$  satisfying this extra requirement is not a problem    unless $n=4$, in which case  we run into  fatal difficulties. But if $n=4$, then in the situation considered here $\delta_\mathscr{C} $ has a maximum at the origin, where $\delta_\mathscr{C}(\vv{0})=2=n-2$, and $\delta_\mathscr{C}(\vv{1})=0$, so $\mathscr{C}$ is a special cell class, hence already handled in the beginning of this section.
\qed
\end{potC}

\section{Acknowledgements}

The author wishes to express his deepest gratitude  to Paul Monsky, for his assistance in the preparation of this paper, for his valuable comments and support, and in particular   for the suggestion of the approach used in Section \ref{bounds}. 

 \appendix
  
 \section{A continuity property of Hilbert--Kunz multiplicities}\label{appendix}
 
 Let $(R,\mm)$ be a Noetherian  local domain of characteristic $p$ and dimension $a\ge 1$, and let $J=\ideal{x_1,\ldots,x_s}$ be an $\mm$-primary ideal of $R$, where $x_1\cdots x_s\ne 0$. 
 
 \begin{defnapp}$J(k)$ is the ideal $\ideal{x_1^k,\ldots,x_s^k}$ of $R$, and $\mu(k)$ is the Hilbert--Kunz multiplicity of $R$ with respect to $J(k)$. 
 \end{defnapp}
 
 It follows immediately  from the definition of the Hilbert--Kunz multiplicity  that $\mu(pk)=p^a\cdot \mu(k)$, so we can extend $\mu$ to a function   $\mu: \mathscr{I}\to \mathbb{R}$, defining
 \[\mu\left(\frac{k}{q}\right)=\frac{\mu(k)}{q^a}.\]
 We shall prove  the following:
 
 \begin{thmapp}\label{mu is lipschitz}
 $\mu$ is a Lipschitz function. In particular, $\mu$ extends uniquely  to a continuous function $\mu^*:[0,1]\to \mathbb{R}$. 
 \end{thmapp} 
 
 We start with a couple of estimates. 
 
\begin{lemapp}\label{estimate 1} $\length_R (J(k-1)/J(k)) =O(k^{a-1})$.
\end{lemapp}

\begin{pf}
Since $J(k-1)/J(k)$ is annihilated by $x_1\cdots x_s$, it is a module over   $S:=R/\ideal{x_1\cdots x_s}$. Let $I$ and $I(k)$ be the extensions of $J$ and $J(k)$ in $S$. Then  $\length_R (J(k-1)/J(k))=\length_S (J(k-1)/J(k))\le \length_S (S/I(k))$, and it suffices to show that this last length is $O(k^{a-1})$. But $I^{sk}\subseteq I(k)$, so $\length_S(S/I(k))\le \length_S(S/I^{sk})$, and the latter is a polynomial in $k$ of degree $a-1$ for $k\gg 0$, since $\dim S=a-1$.
\qed
\end{pf}

\begin{lemapp}\label{estimate 2}
$\mu(k)-\mu(k-1) =O(k^{a-1})$.
\end{lemapp}

\begin{pf} This follows from Lemma \ref{estimate 1} and Lemma 4.2 of \cite{watanabe}, which says that $\mu(k)-\mu(k-1)\le (\text{constant}) \cdot \length_R(J(k-1)/J(k))$, where the constant is the Hilbert--Kunz multiplicity of $R$ with respect to  its maximal ideal $\mm$.
\qed\end{pf}

 The Lipschitz property for $\mu$  follows easily from Lemma \ref{estimate 2}.
 
 \begin{potLipschitz} By Lemma \ref{estimate 2}, there is a constant $M$ such that $\mu(k)-\mu(k-1) \le M k^{a-1}$, for all $k>0$.
 Let $j/q\le k/q$ be two elements of $\mathscr{I}$.  Then 
  $0\le \mu(k)-\mu(j)\le Mk^{a-1}(k-j),$ and 
 dividing by $q^a$ we find   \[0\le \mu\left(\frac{k}{q}\right)-\mu\left(\frac{j}{q}\right)\le M\cdot\left(\frac{k}{q}-\frac{j}{q}\right).\]
\qed \end{potLipschitz}
 
 \begin{remapp}
 With minor modifications in this argument one could prove the following generalization. Let $J(\vv{k})=J(k_1,\ldots,k_s)=\ideal{x_1^{k_1},\ldots,x_s^{k_s}}$ and let $\mu(\vv{k})$ be the Hilbert--Kunz multiplicity of $R$ with respect to $J(\vv{k})$. Then the function 
 \begin{align*}
\mu:  \mathscr{I}^s &\longrightarrow \mathbb{R}\\
 \frac{\vv{k}}{q} &\longmapsto \frac{\mu(\vv{k})}{q^a}
 \end{align*}
 is Lipschitz, and hence can be extended to a continuous function $\mu^*:[0,1]^s\to \mathbb{R}$.
 \end{remapp}

\end{document}